\documentclass[notitlepage, longbibliography, nofootinbib]{revtex4-2}%
\usepackage[hidelinks]{hyperref}
\usepackage{float}
\usepackage{diagbox}
\usepackage{placeins}
\usepackage{pbox}
\usepackage{graphicx}
\usepackage{amsmath}
\usepackage{amssymb}
\usepackage{verbatim}
\usepackage{float}
\usepackage{amsfonts}
\usepackage{tikz}
\usepackage{algorithm}
\usepackage{algpseudocode}%
\providecommand{\U}[1]{\protect\rule{.1in}{.1in}}
\makeatletter

\renewcommand*{\p@subsection}{}

\renewcommand*{\p@subsubsection}{}
\makeatother
\newtheorem{theorem}{Theorem}[section]

\newtheorem{lem}[theorem]{Lemma}

\newtheorem{thm}[theorem]{Theorem}
\newtheorem{conj}[theorem]{Conjecture}

\newtheorem{openq}[theorem]{Open question}
\newtheorem{problem}[theorem]{Problem}

\begin{document}
\title{Attainable bounds for algebraic connectivity and maximally-connected regular graphs}
\author{Geoffrey Exoo}
\email{geoffrey.exoo@gmail.com}
\author{Theodore Kolokolnikov}
\email{tkolokol@gmail.com, corresponding author}
\author{Jeanette Janssen}
\email{Jeannette.Janssen@dal.ca}
\author{Timothy Salamon}
\email{salamontimothy@gmail.com}
\affiliation{$^*$ Department of Computer Science, Indiana State University,
Terre Haute, IN 47809, USA}
\affiliation{$^\dagger {}^\ddagger {}^\S$ Department of Mathematics and Statistics,
Dalhousie University Halifax,
Nova Scotia, B3H3J5, Canada}


\begin{abstract}We derive attainable upper bounds on the algebraic
connectivity (spectral gap) of a regular graph in terms of its diameter and
girth. This bound agrees with the well-known Alon-Boppana-Friedman bound for
graphs of even diameter, but is an improvement for graphs of odd diameter. For
the girth bound, we show that only Moore graphs can attain it, and these only
exist for very few possible girths. For diameter bound, we use a combination
of stochastic algorithms and exhaustive search to find graphs which attain it.
For 3-regular graphs, we find attainable graphs for all diameters $D$ up to
and including $D=9$ (the case of $D=10$ is open). These graphs are extremely
rare and also have high girth; for example we found exactly 45 distinct cubic graphs
on 44 vertices attaining the upper bound when $D=7$; all have girth 8 (out of a total of about $10^{20}$ cubic graphs on 44 vertices, including 266362 having girth 8). We also exhibit families of $d$-regular graphs
attaining upper bounds with $D=3$ and $4$, and with $g=6.$ Several conjectures
are proposed.
\end{abstract}
\maketitle

\section{Introduction}

\tikzstyle{fred}=[draw=black, fill=yellow, shape=circle, minimum height=1.4mm,
inner sep=1, text=black]

The Algebraic Connectivity of a graph (AC; also called the spectral gap) is an
important measure of how well information propagates through the graph
\cite{fiedler1973algebraic, de2007old}, and corresponds to the second
eigenvalue of the graph Laplacian matrix. The higher the AC, the better the
graph is at diffusing information
\cite{olfati2004consensus,olfati2005ultrafast, ghosh2006growing}. Graphs with
high algebraic connectivity are related to expander graphs, and are important
in many applications \cite{hoory2006expander, lubotzky2012expander}. This
paper is motivated by the following question.

\begin{problem}
\label{openq1}\emph{Find graphs which have the highest possible AC among all
regular graphs of a given degree d and diameter }$D,$\emph{ girth }$g,$\emph{
or number of nodes }$n.$
\end{problem}

There are numerous works addressing this and related questions. A well known
upper bound for AC in terms of diameter is the so-called Alon-Boppana-Friedman
bound -- see \cite{friedman1993some, alon1986eigenvalues, nilli2004tight,
hoory2006expander}. Papers \cite{ogiwara2015maximizing, shahbaz2023algebraic,
mosk2008maximum, kolokolnikov2015maximizing} consider the question of
maximizing AC\ over some families of possible graphs with a fixed number of
vertices and edges. In paper \cite{wang2008algebraic} AC is maximized subject
to a fixed diameter and number of edges. Papers
\cite{ghosh2006growing,kim2009bisection,li2018maximizing} propose
graph-growing algorithms to generate a high-AC\ graph on a large network with
a given number of edges and vertices. Paper \cite{kolokolnikov2021better}
considers maximizing AC for several families of random regular graphs. A
complimentary question is explored in \cite{cioaba2016maximizing,
cioabua2020eigenvalues}, which asks what is the largest $n$ for a given AC and
degree $d$.

In this work, we explore question \ref{openq1} for a fixed girth or diameter.
For even diameter, \cite{friedman1993some} gives a tight upper bound which --
as we will see in \S\ref{sec:D} -- is attained in many cases. In theorem
\ref{thm:max} below, we will also derive the analogous tight bound for odd
diameter, as well as for odd and even girths.

\begin{thm}
\label{thm:max}Suppose that a $d-$regular graph has girth $g$ and diameter $D$
and let $AC$ be its algebraic connectivity. Then%
\begin{equation}
AC \leq d-2(d-1)^{1/2}\cos\theta\label{ACmax}%
\end{equation}

where $\theta$ can be two of the following values.

\begin{itemize}
\item If $D$ is even with $D=2K,$ then $\theta$ is the smallest positive root
of
\begin{equation}
\tan(\theta K)=-\frac{d}{d-2}\tan\theta. \label{Deven}%
\end{equation}

\item If $D$ is odd with $D=2K-1,$ then $\theta$ be the smallest positive root
of
\begin{equation}
\tan(\theta K)=-\frac{\left(  2\sqrt{d-1}\cos\theta+d\right)  \sin\theta
}{\sqrt{d-1}\left(  d-2\cos^{2}\theta\right)  +\left(  d-2\right)  \cos\theta
}. \label{Dodd}%
\end{equation}

\item If $g$ is even with $g=2K,$ then $\theta=\pi/K.$

\item If $g$ is odd with $g=2K+1,$ then $\theta$ be the smallest root of
\begin{equation}
\tan\left(  \theta K\right)  =-\frac{\sin\theta}{\left(  d-1\right)
^{-1/2}+\cos\theta}. \label{godd}%
\end{equation}

\end{itemize}
\end{thm}

\begin{figure}[tb]
\includegraphics[width=1\textwidth]{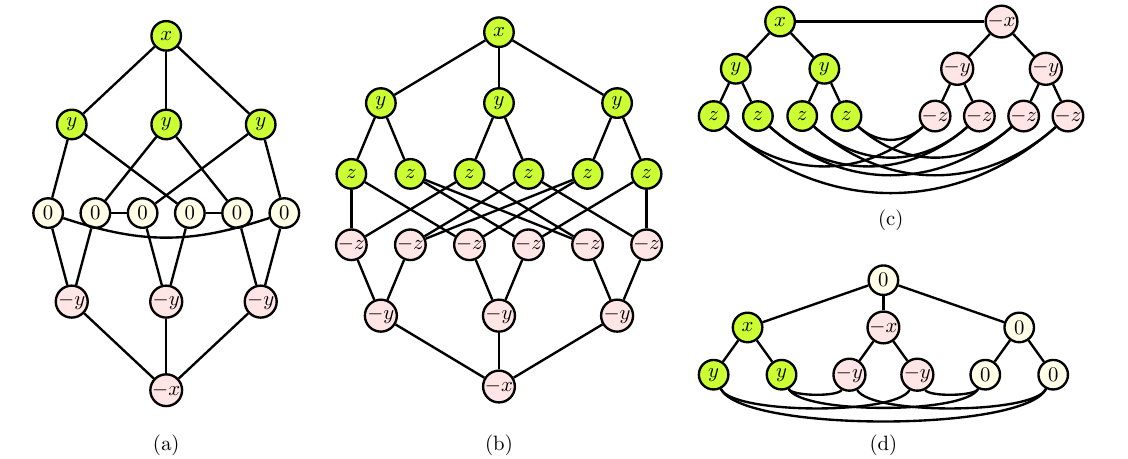}
\caption{Some maximal cubic graphs and the form of the eigenvector
corresponding to the second eigenvalue $\lambda_{2}$ of Laplacian. (a)
Diameter-maximal graph ``Crossing number 3H'' \cite{pegg2009crossing} with
D=4, AC=1.2679. (b) Diameter-maximal graph \textquotedblleft Cubic integral
G10\textquotedblright\ \cite{bussemaker1976there} with D=5, AC=1. (c)
Girth-maximal graph with $g=6$, which is the Heawood graph having AC=1.585.
(d) Girth-maximal graph with $g=5$, which is the Petersen graph having AC=2. }%
\label{fig:max}%
\end{figure}

\begin{table}[tb]%
\[%
\begin{tabular}
[c]{|c|ccccccccc|}\hline
\multicolumn{10}{|c|}{Upper bound for AC in terms of girth}\\\hline
\backslashbox{$g$}{$d$} & 3 & 4 & 5 & 6 & 7 & 8 & 9 & 10 & 11\\\hline
3 & \textbf{4.0000} & \textbf{5.0000} & \textbf{6.0000} & \textbf{7.0000} &
\textbf{8.0000} & \textbf{9.0000} & \textbf{10.000} & \textbf{11.000} &
\textbf{12.000}\\
4 & \textbf{3.0000} & \textbf{4.0000} & \textbf{5.0000} & \textbf{6.0000} &
\textbf{7.0000} & \textbf{8.0000} & \textbf{9.0000} & \textbf{10.000} &
\textbf{11.000}\\
5 & \textbf{2.0000} & \emph{2.6972} & \emph{3.4384} & \emph{4.2087} &
\emph{5.0000} & \emph{5.8074} & \emph{6.6277} & \emph{7.4586} & \emph{8.2984}%
\\
6 & \textbf{1.5858} & \textbf{2.2679} & \textbf{3.0000} & \textbf{3.7639} &
\emph{4.5505} & \textbf{5.3542} & \textbf{6.1716} & \textbf{7.0000} & 7.8377\\
7 & \emph{1.1864} & \emph{1.7466} & \emph{2.3738} & \emph{3.0443} &
\emph{3.7458} & \emph{4.4709} & \emph{5.2147} & \emph{5.9739} & \emph{6.7460}%
\\
8 & \textbf{1.0000} & 1.5505 & 2.1716 & 2.8377 & 3.5359 & 4.2583 & 5.0000 &
5.7574 & 6.5279\\
9 & \emph{0.8088} & \emph{1.3004} & \emph{1.8706} & \emph{2.4913} &
\emph{3.1481} & \emph{3.8322} & \emph{4.5380} & \emph{5.2616} & \emph{6.0000}%
\\
10 & \emph{0.7118} & \emph{1.1975} & \emph{1.7639} & \emph{2.3820} &
\emph{3.0366} & \emph{3.7191} & \emph{4.4235} & \emph{5.1459} & \emph{5.8833}%
\\
11 & \emph{0.6069} & \emph{1.0600} & \emph{1.5983} & \emph{2.1912} &
\emph{2.8229} & \emph{3.4840} & \emph{4.1685} & \emph{4.8721} & \emph{5.5916}%
\\
12 & \textbf{0.5505} & 1.0000 & 1.5359 & 2.1270 & 2.7574 & 3.4174 & 4.1010 &
4.8038 & 5.5228\\
13 & \emph{0.4872} & \emph{0.9168} & \emph{1.4356} & \emph{2.0114} &
\emph{2.6277} & \emph{3.2748} & \emph{3.9462} & \emph{4.6376} & \emph{5.3456}%
\\\hline
\end{tabular}
\]%
\[%
\begin{tabular}
[c]{|c|ccccccccc|}\hline
\multicolumn{10}{|c|}{Upper bound for AC in terms of diameter}\\\hline
\backslashbox{$D$}{$d$} & 3 & 4 & 5 & 6 & 7 & 8 & 9 & 10 & 11\\\hline
3 & \textbf{2.0000} & \textbf{3.0000} & \textbf{4.0000} & \textbf{5.0000} &
\textbf{6.0000} & \textbf{7.0000} & \textbf{8.0000} & \textbf{9.0000} &
\textbf{10.000}\\
4 & \textbf{1.2679} & \textbf{2.0000} & \textbf{2.7639} & \textbf{3.5505} &
\textbf{4.3542} & \textbf{5.1716} & \textbf{6.0000} & \textbf{6.8377} &
7.6834\\
5 & \textbf{1.0000} & \textbf{1.6972} & 2.4384 & 3.2087 & 4.0000 & 4.8074 &
5.6277 & 6.4586 & 7.2984\\
6 & \textbf{0.7639} & 1.3542 & 2.0000 & 2.6834 & 3.3944 & 4.1270 & 4.8769 &
5.6411 & 6.4174\\
7 & \textbf{0.6571} & 1.2266 & 1.8587 & 2.5321 & 3.2356 & 3.9621 & 4.7070 &
5.4671 & 6.2398\\
8 & \textbf{0.5505} & 1.0665 & 1.6508 & 2.2810 & 2.9446 & 3.6340 & 4.3440 &
5.0709 & 5.8123\\
9 & \textbf{0.4965} & 1.0000 & 1.5762 & 2.2006 & 2.8597 & 3.5456 & 4.2527 &
4.9772 & 5.7164\\
10 & 0.4384 & 0.9111 & 1.4601 & 2.0598 & 2.6964 & 3.3613 & 4.0487 & 4.7546 &
5.4762\\
11 & 0.4069 & 0.8717 & 1.4156 & 2.0118 & 2.6456 & 3.3083 & 3.9939 &
4.6984 & 5.4187\\
12 & 0.3714 & 0.8167 & 1.3436 & 1.9245 & 2.5444 & 3.1941 & 3.8677 & 4.5607 &
5.2701\\
13 & 0.3512 & 0.7912 & 1.3148 & 1.8934 & 2.5115 & 3.1599 & 3.8323 & 4.5243 &
5.2329\\\hline
\end{tabular}
\]
\caption{Upper bounds for AC in terms of girth and diameter. Known attainable
bounds are in bold. Known unattainable are in italics. The rest are unknown.}%
\label{table:max}%
\end{table}The derivation of this theorem is given in \S\ref{sec:proofs}. For even diameter, the bound (\ref{Deven})\ was already obtained in
\cite{friedman1993some} (see Proposition 3.2 there; equation (\ref{Deven})\ is
equivalent to formula for $\theta$ in Corollary 3.6 of \cite{friedman1993some}%
). For odd diameter, formula (\ref{Dodd}) is an improvement to the current
literature as far as we are aware. The tight upper bound $\theta=\pi/K$ for
even girth was previously derived in \cite{kolokolnikov2015maximizing} in the
context of cubic graphs. The bound for odd girth was given in
\cite{salamon2022algebraic}.

For $g,G\leq6$ and arbitrary $d,$ we have the following explicit formulas
upper bounds:%
\[%
\begin{tabular}
[c]{|l|l|}%
\hline
$D$ & $AC$ upper bound\\\hline\hline
3 & $d-1$\\\hline
4 & $d-\sqrt{d}$\\\hline
5 & $d-\frac{1}{2}-\sqrt{d-\frac{3}{4}}$\\\hline
6 & $d-\sqrt{2d-1}$\\
\hline
\end{tabular}
~~~~~~~~~~~~
\begin{tabular}
[c]{|c|l|}%
\hline
$\ \ g\ \ $ & $AC$ upper bound\\\hline\hline
3 & $d+1$\\\hline
4 & $d$\\\hline
5 & $d+\frac{1}{2}-\sqrt{d-\frac{3}{4}}$\\\hline
6 & $d-\sqrt{d-1}$\\
\hline
\end{tabular}
\]

It is easy to see that in all four cases of Theorem \ref{thm:max}, $\frac{\pi
}{2K}<\theta<\frac{\pi}{K}$, and $\theta\rightarrow\frac{\pi}{K}$ as
$K\rightarrow\infty.$ The bound $\theta=\pi/K$ was also derived by Nilli
\cite{nilli2004tight}, and both asymptote to the the Alon-Boppana estimate
\cite{friedman1991second, friedman1993some} of $AC\sim d-2\sqrt{d-1}$ for
random regular graphs as $D\rightarrow\infty$.

Table \ref{table:max} lists the upper bounds for small $d,g,D$, and summarizes
known attainable bounds. Many of these bounds are attainable, particularly
with respect to diameter. Let us define a graph to be
\textbf{diameter-maximal} or \textbf{girth-maximal} if it attains diameter
(resp. girth) bound of theorem \ref{thm:max} (we will abbreviate it to
\textquotedblleft maximal\textquotedblright\ when prefix is understood from
the context). Sections \ref{sec:D} and \ref{sec:girth} of this paper are dedicated
to a search for maximal graphs.

To illustrate where these bounds come from, consider the graph in figure
\ref{fig:max}(a). An eigenvector assigns a number to each vertex. Here, we
take a specific eigenvector whose entries have values of $\pm x,\pm y,0$
associated to vertices as shown in the figure. This choice guarantees that the
sum of the entries in the eigenvector is zero, and as such it is orthogonal to
the $\left[  1,1,\ldots,1\right]  ^{T}$ eigenvector corresponding to the zero
eigenvalue of the Laplacian. Moreover both top and bottom vertices yield
$\lambda x=3x-3y;$ whereas vertices at 2nd and 4th rows both read $\lambda
y=3y-x.$ Correspondingly, the Laplacian matrix includes the spectrum of the
matrix $\left[
\begin{array}
[c]{rr}%
3 & -3\\
-1 & 3
\end{array}
\right]  .$ The smallest eigenvalue of this matrix is $1.2679,$ which is the
upper bound for the AC of cubic graphs of diameter 4; a bound that is in fact
attained by this graph.

Similarly, vertex assignment in (b), (c) and (d) yields matrices
\[
\left[
\begin{array}
[c]{rrr}%
3 & -3 & 0\\
-1 & 3 & -2\\
0 & -1 & 5
\end{array}
\right]  ,\ \left[
\begin{array}
[c]{rrr}%
4 & -2 & 0\\
-1 & 3 & -2\\
0 & -1 & 5
\end{array}
\right]  ,\ \text{and\ }\left[
\begin{array}
[c]{rr}%
3 & -2\\
-1 & 5
\end{array}
\right]
\]
respectively, whose smallest eigenvalues (and upper bounds for AC)\ are
$1,\ 1.5858,$ and$\ 2,$ respectively. In fact, these numbers are exactly the
AC\ for these graphs.

Naturally most graphs don't have such a nice structure. Nonetheless we show in
\S\ref{sec:proofs} that these bounds hold for all graphs. For maximal
graphs with odd diameter or odd/even girth, the proof also yields the exact
number of vertices needed (the case of even diameter is more complicated). For
example, figure \ref{fig:max}(b)\ shows a maximal diameter-5 cubic graph with
$n=20$ vertices. More generally, when $D=2K-1$ is odd, any maximal graph must
contain exactly two Moore trees with $K$ levels. Girth-maximal graphs must
consist entirely of the corresponding Moore trees (rooted at a vertex for odd
girth, or at an edge for an even girth) along with additional edges joining
the leafs. Correspondingly, maximal-girth graphs can only be attained by Moore
graphs. This puts a severe restriction on the possible girth-maximal graphs.
We summarize this as follows.

\begin{theorem}
\label{thm:Dodd}Suppose that a diameter-maximal graph has an odd diameter,
$D=2K-1.$ Then it is bipartite, and consists of two disjoint Bethe trees of
$K$ levels with an additional $d\left(  d-1\right)  ^{K-1}$ edges joining
their leafs. It has exactly
\begin{equation}
n=2\left(  1+d\sum_{j=0}^{K-2}\left(  d-1\right)  ^{j}\right)  =2\left(
\frac{d\left(  d-1\right)  ^{K-1}-2}{d-2}\right)
\end{equation}
vertices.
\end{theorem}

\begin{theorem}
\label{thm:girth}A girth-maximal graph must be a Moore graph, that is, a graph
that attains a Moore bound for the number of vertices.
\end{theorem}

The situation is more complicated for even diameter $D=2K.$ In this case, a
maximal graph consists of two Bethe trees with $K$ levels, plus a middle layer
that is not part of either tree (see figure \ref{fig:max}(a)). Each Bethe
tree has $d-1$ edges going from each of its $d\left(  d-1\right)  ^{K-2}$
leafs to the middle layer, for a total of $2d\left(  d-1\right)  ^{K-1}$ edges
going into the middle layer. This requires a minimum of $2\left(  d-1\right)
^{K-1}$ of vertices in the middle layer. We summarize:

\begin{theorem}
\label{thm:Deven}Suppose that a diameter-maximal graph has an even diameter,
$D=2K.$ Then it consists of two disjoint Bethe trees of $K$ levels, plus a
center of at least $\left(  d-1\right)  ^{K-1}$ vertices that is not part of
either trees. The total number of vertices $n$ satisfies%
\[
n\geq2\left(  1+d\sum_{j=0}^{K-2}\left(  d-1\right)  ^{j}\right)  +2\left(
d-1\right)  ^{K-1}=4\left(  \frac{\left(  d-1\right)  ^{K}-1}{d-2}\right)  .
\]

\end{theorem}

This lower bound on $n$ is sometimes attained and sometimes not. For
example,when $d=3,D=4$,\ the lower bound $n=12$ is \emph{not} attained (the smallest
such maximal graph has 14 vertices, see \S \ref{sec:D34}). On the other hand,
the lower bound is $n=16$ when $d=4,D=4,$ and it\emph{ is} attained by six
distinct graphs (see \S \ref{sec:D44}).

Unlike odd diameter, we do not have a tight upper bound on $n$ for maximal
graphs with even diameter. Asides from the obvious Moore bound, a nontrivial
upper bound is given in \cite{cioabua2020eigenvalues}. For example consider
the case of $d=3,D=4$, for which Theorem \ref{thm:max} gives $\max
AC\approx1.2679.$ Applying Theorem 8 from \cite{cioabua2020eigenvalues} with
$r=3,t=5$ and $c=4.1$ yields $M\approx27.8$ and therefore an upper bound of
$n\leq26,$ which is the best possible bound that can be obtained from Theorem
8 of \cite{cioabua2020eigenvalues} for this case. On the other hand, the
biggest maximal graph that we were able to find (and the conjectured upper
bound)\ has $n=18$ vertices (see \S\ref{sec:D34}).

The outline of the paper is as follows. Theorems
\ref{thm:max} to \ref{thm:Deven} are derived in \S \ref{sec:proofs}. In \S \ref{sec:D} we classify maximal
graphs for specific values of $D$ and $d.$ In \S \ref{sec:girth} we discuss
maximizing AC for a given girth. In \S \ref{sec:families} we present several
infinite families of infinite graphs for $D=3,4$ and $g=6$.  Some open questions and
conjectures are posed in \S\ref{sec:discuss}.

\section{Proofs of Theorems \ref{thm:max}-\ref{thm:Deven}} \label{sec:proofs}

In this section we derive the upper bounds in Theorem \ref{thm:max}. The bounds on graph order in Theorems \ref{thm:Dodd}-\ref{thm:Deven} will follow from the proof of Theorem \ref{thm:max}.

The proof
of Theorem \ref{thm:max} follows closely the arguments presented in
\cite{nilli2004tight, friedman1993some}, with an additional argument to get a tighter bound in the case of odd diameter.

We start with the following lemma.

\begin{lem}
Consider a $K\times K$ tri-diagonal matrix%
\begin{equation}
M=\left[
\begin{array}
[c]{ccccc}%
a & -b &  &  & \\
-1 & d & -\left(  d-1\right)   &  & \\
& -1 & \ddots & \ddots & \\
&  & \ddots & d & -\left(  d-1\right)  \\
&  &  & -1 & c
\end{array}
\right]  .
\end{equation}
Its eigenvalues are given by%
\begin{equation}
\lambda=d-(d-1)^{1/2}2\cos\theta\label{lam}%
\end{equation}
where $\theta\neq0$ satisfies the following equation, depending on values of
$a,b,c.$

\begin{itemize}
\item[(a)] If $a=d,b=d$ and $c=d$ then%
\begin{equation}
\tan\left(  \theta K\right)  =-\frac{d}{d-2}\tan\theta\label{theta-even}%
\end{equation}

\item[(b)] If $a=d,b=d$ and $c=2d-1$ then%
\begin{equation}
\tan(\theta K)=\frac{\left(  2\sqrt{d-1}\cos\theta+d\right)  \sin\theta}%
{\sqrt{d-1}\left(  2\cos^{2}\theta-d\right)  +\left(  2-d\right)  \cos\theta}.
\label{theta-odd}%
\end{equation}

\item[(c)] If $a=d+1,b=d-1,\ $and $c=2d+1$ then
\begin{equation}
\sin\left(  K\theta\right)  =0\label{theta-girth-even}%
\end{equation}

\item[(d)] If $a=d,b=d-1,\ $and $c=d+1$ then
\begin{equation}
\tan\left(  \theta K\right)  =-\frac{\sin\theta}{\left(  d-1\right)
^{-1/2}+\cos\theta}\label{theta-girth-odd}%
\end{equation}

\end{itemize}

For the case (a), the corresponding eigenvector has entries%
\begin{equation}
v_{j+1}=\left(  d-1\right)  ^{-j/2}\sin\left(  \theta(j-K)\right)
,\ \ j=0\ldots K-1.\label{v-casea}%
\end{equation}
For the case (b), the corresponding eigenvector has entries%
\begin{equation}
v_{j+1}=\left(  d-1\right)  ^{-j/2}\left[  \sin\left(  \theta(j-K)\right)
+(d-1)^{1/2}\sin\left(  \theta(j+1)\right)  \right]  ,\ \ j=0\ldots K-1.
\end{equation}
In both cases, $v_{j}$ decreases with $j.$

In all cases the smallest eigenvalue satisfies $\pi/2<\theta K\leq\pi.$
\end{lem}

\textbf{Proof.}

We consider the following anzatz for the eigenvector,%
\begin{equation}
v_{j+1}=\left(  d-1\right)  ^{-j/2}\left[  Ae^{i\theta j}+Be^{-i\theta
j}\right]  ,\ j=0\ldots K-1.\label{304}%
\end{equation}
Then all rows except the first and last yield the same equation, namely.%
\[
\lambda=d-(d-1)^{1/2}\left(  e^{i\theta}+e^{-i\theta}\right)  .
\]
Consider first the case $a,b=d.$ Then the first row simplifies to
\begin{equation}
A\left\{  e^{i\theta}-\left(  d-1\right)  e^{-i\theta}\right\}  +B\left\{
e^{-i\theta}-\left(  d-1\right)  e^{i\theta}\right\}  =0\label{300}%
\end{equation}
while the last row simplifies to%
\begin{equation}
\frac{c-d}{(d-1)^{1/2}}\left[  Ae^{i\theta\left(  K-1\right)  }+Be^{-i\theta
\left(  K-1\right)  }\right]  +Ae^{i\theta K}+Be^{-i\theta K}=0.\label{302}%
\end{equation}
A solution to (\ref{300})\ is given by%
\begin{equation}
A=e^{-i\theta}-\left(  d-1\right)  e^{i\theta},\ \ \ B=-\bar{A}.\label{301}%
\end{equation}
Upon substituting (\ref{301}) into (\ref{302})\ we obtain%
\[
\sin\left(  \theta\left(  K-1\right)  \right)  -\left(  d-1\right)
\sin\left(  \theta\left(  1+K\right)  \right)  +\frac{c-d}{(d-1)^{1/2}}\left(
\sin\left(  \theta\left(  K-2\right)  \right)  -\left(  d-1\right)
\sin\left(  \theta K\right)  \right)  =0.
\]
Using the sine addition formula then yields%
\begin{equation}
\tan\left(  \theta K\right)  =\frac{\left(  2\left(  c-d\right)  d\cos
\theta+d\sqrt{d-1}\right)  \sin\theta}{\left(  c-d\right)  \left(  2\cos
^{2}\theta-d\right)  d+\left(  2-d\right)  \sqrt{d-1}\cos\theta}.\label{303}%
\end{equation}
Cases (\ref{theta-even}) and\ (\ref{theta-odd}) correspond to special cases of
(\ref{303}).

Substituting (\ref{301})\ into (\ref{304})\ yields the formula for the
eigenvector
\begin{equation}
v_{j+1}=\left(  d-1\right)  ^{-j/2}\left[  \sin\left(  \theta(j-K)\right)
+\frac{c-d}{(d-1)^{1/2}}\sin\left(  \theta(j+1)\right)  \right]
,\ \ j=0\ldots K-1.\label{305}%
\end{equation}
Formulas (a)\ and (b) correspond to special cases of (\ref{305}).

Other cases are derived similarly; we omit the details. $\blacksquare$

\textbf{Even diameter. }We now present the proof of Theorem \ref{thm:max} for the case of the even
diameter. To illustrate the proof, consider the graph such as in figure
\ref{fig:max}(a). It consists of a \textquotedblleft
double-tree\textquotedblright\ structure plus a middle layer. The top and
bottom trees both have $K=D/2$ levels (here, $K=2,D=4$). The vertices in each
level of the top tree are assigned the same value, and the opposite value is
assigned to layers of the bottom tree. Vertices in the middle layer are
assigned a value of zero. For our symmetric example, this guarantees that the
entries of the resulting vector sum to zero, so that it is perpendicular to
the eigenvector $\left(  1,1,\ldots1\right)  ,$ and the corresponding
eigenvalue therefore bounds AC.

Now most graphs are not so symmetric as the one shown in figure \ref{fig:max}.
Nonetheless we can still use the same assignment as a test vector. The idea is
to consider the top and bottom tree separately, then combine them together to
get a bound for AC. To this end, consider a graph $T$ where all the vertices
are a distance at most $K$ away from vertex 1, and consider the following
eigenvalue problem on such a graph:%
\begin{equation}
\left\{
\begin{array}
[c]{c}%
\lambda x_{k}=\sum\limits_{j\sim k}(x_{k}-x_{j}),\ \text{if }%
\operatorname{dist}(k,1)<K\\
\lambda x_{k}=(d-\deg(k))x_{k}+\sum\limits_{j\sim k}(x_{k}-x_{j})\text{,
}\operatorname{dist}(k,1)=K
\end{array}
\right.  ,\label{1235}%
\end{equation}
\textbf{ }where $j\sim k$ if vertices $j,k$ are neighbours, $\deg(k)$ denotes
the degree of vertex $k$, and $\operatorname{dist}(k,1)$ is the distance of
vertex $k$ from vertex 1. Let $T_{K,d}$ be a Bethe tree with $K$ levels of
degree $d$, tree as illustrated below:%
\[%
\begin{array}
[c]{c}%
T_{K,d}\text{ with }K=4,d=3\\
\text{\includegraphics[width=0.5\textwidth]{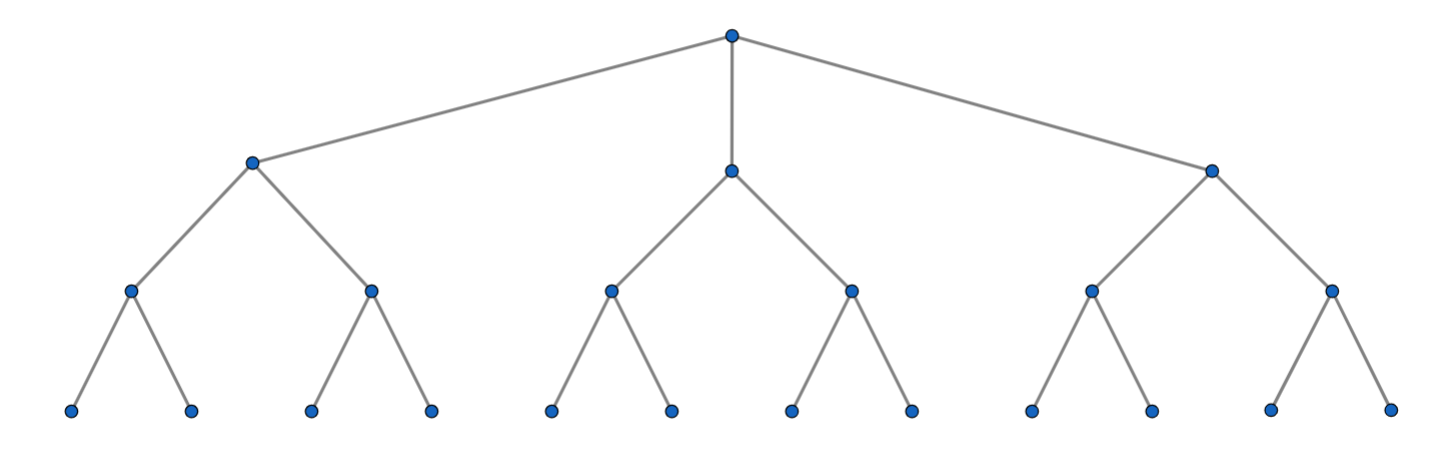}}%
\end{array}
\]
In other words, a tree having $K$ levels where each non-leaf nodes has degree
$d.$ We have the following lemma.

\begin{lem}
\label{lem:even}Suppose that $T$ is a where all the vertices are a distance at
most $K$ away from vertex 1, and where each vertex has degree at most $d,$ and
let $\lambda_{\min}(T)$ be the smallest eigenvalue of $T.$ Then $\lambda
_{\min}(T)\leq\lambda_{\min}(T_{K,d})$ where $T_{K,d}$ is the Moore tree of
level $K$ and degree $d.$ Moreover, the equality only happens when
$T=T_{K,d}.$ Explicitly, $\lambda_{\min}(T_{K,d})$ is given by\ (\ref{lam})
where $\theta>0$ is the smallest solution of (\ref{theta-even}).
\end{lem}

\textbf{Proof. }The argument we present here is essentially the same as that
in \cite{nilli2004tight}. First, we compute $\lambda=\lambda_{\min}(T_{K,d}).$
The corresponding eigenvector is obtained by assigning the same value
$v_{k+1}$ to all vertices distance $k$ from the root, $k=0\ldots K-1.$ It
therefore satisfies
\begin{align}
\lambda v_{1} &  =d\left(  v_{1}-v_{2}\right)  \label{325}\\
\lambda v_{2} &  =dv_{2}-v_{1}-\left(  d-1\right)  v_{3}\nonumber\\
&  \vdots\\
\lambda v_{K} &  =dv_{K}-v_{k-1}\nonumber
\end{align}
By Lemma \ref{lem:even}, we find that $\lambda=\lambda_{\min}(T_{K,d})$ with
the corresponding eigenvector $v_{j}$ given by \ref{v-casea}.

Next, we take a test vector by assigning $v_{j}$ to all vertices at the level
$j$ $\left(  j=1\ldots K\right)  $ of $T.$ The associated Rayleigh quotient
is\
\begin{equation}
R=\frac{\sum_{j=1}^{K-1}e_{j}\left(  v_{j}-v_{j+1}\right)  ^{2}+e_{K}v_{K}%
^{2}}{\sum_{j=1}^{K}n_{j}v_{j}^{2}}%
\end{equation}
where $e_{j-1}$ is the number of edges from level $j-1$ to level $j$; and
$e_{K}=dn_{K}-e_{K-1}.$ Since each vertex has degree at most $d,$ we also
have:
\[
e_{1}\leq d;\ \ \ \ e_{j}+e_{j-1}\leq dn_{j},\ \ j=2\ldots K
\]
so that
\begin{align*}
\sum_{j=1}^{K-1}e_{j}\left(  v_{j}-v_{j+1}\right)  ^{2}+e_{K}v_{K}^{2} &
=e_{1}v_{1}^{2}+\sum_{j=2}^{K}\left(  e_{j}+e_{j-1}\right)  v_{j}^{2}%
-2\sum_{1}^{K-1}e_{j}v_{j}v_{j+1}\\
&  \leq d\sum_{j=1}^{K}n_{j}v_{j}^{2}-\sum_{1}^{K-1}2e_{j}v_{j}v_{j+1}.
\end{align*}
Write%
\[
\sum_{1}^{K-1}2e_{j}v_{j}v_{j+1}=e_{1}v_{1}v_{2}+e_{K}v_{K-1}v_{K}+\sum
_{2}^{K-1}v_{j}\left(  e_{j-1}v_{j-1}+e_{j}v_{j+1}\right)  .
\]
and consider the term $e_{j-1}v_{j-1}+e_{j}v_{j+1}.$ Note that $e_{j-1}%
+e_{j}=n_{j}d-e_{j\rightarrow j}$ where $e_{j\rightarrow j}$ is the number of
edges that are internal to level $j$ of the tree. In addition, $e_{j-1}\geq
n_{j}$ (with equality when each vertex has only one parent). Moreover, $v_{j}$
is decreasing in $j.$ It follows that the minimum possible value of
$e_{j-1}v_{j-1}+e_{j}v_{j+1}$ corresponds to the case when $e_{j-1}%
=n_{j},\ \ e_{j}=(d-1)n_{j},$ $e_{j\rightarrow j}=0.$ This yields%
\[
e_{j-1}v_{j-1}+e_{j}v_{j+1}\geq n_{j}v_{j-1}+(d-1)n_{j}v_{j+1}=n_{j}%
(d-\lambda)v_{j},\ \ \ j=2\ldots K-1.
\]
Also, $e_{1}=dn_{1}$ so that%
\[
e_{1}v_{2}=dn_{1}v_{2}=n_{1}\left(  d-\lambda\right)  v_{1}.
\]
Finally$\ e_{K}\geq n_{K}$ so that%
\[
e_{K}v_{K-1}\geq n_{K}v_{K-1}=n_{K}\left(  d-\lambda\right)  v_{K}.
\]
Combining above we obtain%

\[
R\leq d-\frac{\sum_{j=1}^{K}n_{j}(d-\lambda)x_{j}^{2}}{\sum_{j=1}^{K}%
n_{j}x_{j}^{2}}=\lambda.
\]

$\blacksquare$

We are now ready to complete the proof of Theorem \ref{thm:max} part\ (a).
Take two vertices $r,\tilde{r}$ separated by a distance $D=2K.$ All vertices
that are at a distance $j$ from $r,$ with $j=0\ldots K-1$, are assigned a
value of $v_{j+1}$, where $v_{j}$ is the eigenvector for the full Bethe tree
as in Lemma \ref{lem:even}. All vertices that are at a distance $j$ from
$\tilde{r},$ with $j=0\ldots K-1$, are assigned a value of $-\alpha v_{j+1}.$
All other vertices are assigned a value of zero. The constant $\alpha$ is
chosen so that the sum of all the assigned values is zero. Then the Rayleigh
quotient\ is given by%
\[
R=\frac{\sum_{j=1}^{K-1}e_{j}\left(  v_{j}-v_{j+1}\right)  ^{2}+e_{K}v_{K}%
^{2}+\alpha^{2}\left(  \sum_{j=1}^{K-1}\tilde{e}_{j}\left(  v_{j}%
-v_{j+1}\right)  ^{2}+\tilde{e}_{K}v_{K}^{2}\right)  }{\sum_{j=1}^{K}%
n_{j}v_{j}^{2}+\alpha^{2}\sum_{j=1}^{K}\tilde{n}_{j}v_{j}^{2}}%
\]
where

\begin{itemize}
\item $n_{j}$ (respectively $\tilde{n}_{j})$ is the number of vertices that
that are assigned weight $v_{j}$ (respectively $\tilde{v}_{j}$), $j=1\ldots
K;$

\item $e_{j}$ (respectively $\tilde{e}_{j})$ is the number of edges between
vertices that have weight $v_{j}$ and $v_{j+1}$ (respectively $\tilde{v}_{j}$
and $\tilde{v}_{j+1}$), $j=1\ldots K-1;$

\item $e_{K}$ (respectively $\tilde{e}_{K})$ is the number of edges between
vertices that have weight $v_{K}$ (respectively $\tilde{v}_{K})$ and zero.
\end{itemize}

We then have,%

\begin{align*}
R &  \leq\min\left(  \frac{\sum_{j=1}^{K-1}e_{j}\left(  v_{j}-v_{j+1}\right)
^{2}+e_{K}v_{K}^{2}}{\sum_{j=1}^{K}n_{j}v_{j}^{2}},\frac{\sum_{j=1}%
^{K-1}\tilde{e}_{j}\left(  v_{j}-v_{j+1}\right)  ^{2}+\tilde{e}_{K}v_{K}^{2}%
}{\sum_{j=1}^{K}\tilde{n}_{j}v_{j}^{2}}\right)  \\
&  \leq\lambda_{\min}(T_{K,d}).
\end{align*}
The first inequality is true since $R$\ is a monotone function of $\alpha^{2}%
$; the second follows from Lemma \ref{lem:even}. This completes the proof of
\ref{thm:max}, part(a).

\textbf{Odd diameter.} The difference between odd and even diameter is
illustrated in \ref{fig:max}(a)\ (even)\ and \ref{fig:max}(b) (odd). The main
difference is that there are edges between the two trees without having the
middle layer. Correspondingly, the problem (\ref{lameven})\ is replaced by the problem%

\begin{equation}
\left\{
\begin{array}
[c]{c}%
\lambda x_{k}=\sum\limits_{j\sim k}(x_{k}-x_{j}),\ \text{if }%
\operatorname{dist}(k,1)<K\\
\lambda x_{k}=2(d-\deg(k))x_{k}+\sum\limits_{j\sim k}(x_{k}-x_{j})\text{,
if}\operatorname{dist}(k,1)=K.
\end{array}
\right.  \label{lameven}%
\end{equation}

The analogue of Lemma \ref{lem:even} still holds with (\ref{325})\ replaced
by
\begin{align}
\lambda v_{1} &  =d\left(  v_{1}-v_{2}\right)  \nonumber\\
\lambda v_{2} &  =dv_{2}-v_{1}-\left(  d-1\right)  v_{3}\nonumber\\
&  \vdots\\
\lambda v_{K-1} &  =dv_{K-1}-v_{K-2}-\left(  d-1\right)  v_{K},\nonumber\\
\lambda v_{K} &  =(2d-1)v_{K}-v_{k-1}\nonumber
\end{align}
As in the even case, take two vertices $r,\tilde{r}$ that are separated by a
distance $D=2K-1.$ All vertices that are at a distance $j$ from $r,$ with
$j=0\ldots K-1$, are assigned a value of $v_{j+1}$. All vertices that are at a
distance $j$ from $\tilde{r},$ with $j=0\ldots K-1$, are assigned a value of
$-\alpha v_{j+1}.$ All other vertices are assigned a value of zero. The
constant $\alpha$ is chosen so that the sum of all the assigned values is
zero. Then the rayleygh quotient then reads,%
\begin{equation}
R\leq\frac{\sum_{j=1}^{K-1}e_{j}\left(  v_{j}-v_{j+1}\right)  ^{2}+e_{K}%
v_{K}^{2}+\alpha^{2}\left(  \sum_{j=1}^{K-1}\tilde{e}_{j}\left(  v_{j}%
-v_{j+1}\right)  ^{2}+\tilde{e}_{K}v_{K}^{2}\right)  +e_{B}\left(
1+\alpha\right)  ^{2}v_{K}^{2}}{\sum_{j=1}^{K}n_{j}v_{j}^{2}+\alpha^{2}%
\sum_{j=1}^{K}\tilde{n}_{j}v_{j}^{2}}.\label{328}%
\end{equation}
Here $n_{j},e_{j}$ are as before, whereas $e_{B}$ is the number of edges
between vertices that have weight $v_{K}$ and $\tilde{v}_{K}$.

Next, note that
\begin{equation}
\left(  1+\alpha\right)  ^{2}\leq2+2\alpha^{2}\label{ineq-alpha}%
\end{equation}
for all $\alpha$, with equality if and only if $\alpha=1$. Replacing $\left(
1+\alpha\right)  ^{2}$ by $2+2\alpha^{2}$ in (\ref{328})\ we therefore obtain%
\[
R\leq\frac{\sum_{j=1}^{K-1}e_{j}\left(  v_{j}-v_{j+1}\right)  ^{2}+\left(
e_{K}+2e_{B}\right)  v_{K}^{2}+\alpha^{2}\left(  \sum_{j=1}^{K-1}\tilde{e}%
_{j}\left(  v_{j}-v_{j+1}\right)  ^{2}+\left(  \tilde{e}_{K}+2e_{B}\right)
v_{K}^{2}\right)  }{\sum_{j=1}^{K}n_{j}v_{j}^{2}+\alpha^{2}\sum_{j=1}%
^{K}\tilde{n}_{j}v_{j}^{2}}.
\]
The resulting expression is monotone in $\alpha^{2}$. It follows that%
\[
R\leq\max\left(  R_{1},\tilde{R}_{1}\right)
\]
where%
\[
R_{1}=\frac{\sum_{j=1}^{K-1}e_{j}\left(  v_{j}-v_{j+1}\right)  ^{2}+\left(
e_{K}+2e_{B}\right)  v_{K}^{2}}{\sum_{j=1}^{K}n_{j}v_{j}^{2}},\ \ \ \tilde
{R}_{1}=\frac{\sum_{j=1}^{K-1}\tilde{e}_{j}\left(  v_{j}-v_{j+1}\right)
^{2}+\left(  \tilde{e}_{K}+2e_{B}\right)  v_{K}^{2}}{\sum_{j=1}^{K}\tilde
{n}_{j}v_{j}^{2}}.
\]
Using the argument identical to Lemma \ref{lem:even}, we find that
$R_{1},\tilde{R}_{1}\leq\lambda$ given by (\ref{lam}, \ref{theta-odd}).

The proof for girth is analogous and is in fact easier. A\ graph of odd girth
$D=2K+1$ has to contain a Bethe tree $T_{K,d}$. This leads to the
corresponding thresholds in Theorem \ref{thm:max}. Similarly a graph of even
girth $D=2K$ has to contain a Moore tree rooted at an edge such as illustrated
in Figure \ref{fig:max}(c) for the case of $K=3$. Detailed proof for even
girth was given in \cite{kolokolnikov2015maximizing}, and in the case of odd
girth, the threshold (\ref{godd}) was worked out in
\cite{salamon2022algebraic}.$\blacksquare$

To show Theorem \ref{thm:Dodd}, one simply traces the inequalities and note
equality is only possible when the two trees coming from $r$ and $\tilde{r}$
are full, and moreover all the edges from the leafs of one tree go to the
edges of the other (corresponding to the extreme value of $\alpha=1$ in
(\ref{ineq-alpha})). Similar arguments show Theorems \ref{thm:Deven} and
\ref{thm:girth}.

\section{Diameter-maximal graphs:\ specific values.\label{sec:D}}

In this section we present graphs that achieve the diameter bounds in Theorem
\ref{thm:max} for small values of $d$ and $D$. We found experimentally that
maximal graphs of diameter $D$ have girth at least $D$ (conjecture
\ref{conj:Dg}).
Subsequently, for graphs of order $n>20$, we
restricted our search to graphs of high girth relative to given order $n$. For
cubic and quartic graphs and sufficiently small $n,$ a complete list of such
graphs is available \cite{coolsaet2023house, mckay1998fast, meringer1999fast}
and in that case we give an exhaustive list of maximal graphs. Where
complete enumeration is impossible (e.g. $D\geq8$ for cubic graphs), we used
a stochastic algorithm (see Appendix A) to search for high-girth graphs. All
of the graphs as well as some programs we used are available for download from
author's website \cite{code}. 

\subsection{Degree 3, Diameter 3\label{sec:D33}}

In this case, the AC bound is $2$, a value achieved only by the graph of the
$3$-cube. The diameter-3 maximal graph is unique for any $d$; see
\S \ref{sec:D3} for the proof.

\subsection{Degree 3, Diameter 4\label{sec:D34}}

\begin{figure}[H]
\centering
\begin{tikzpicture}[scale=0.5,line width=1pt]
\newcommand\orad{4.0}
\newcommand\dtheta{360.0/14.0}
\begin{scope}[rotate=90]
\foreach \theta in {0,1,2,...,13}
{
\begin{scope}[rotate={\theta*\dtheta}]
\draw[fred] (0:\orad) -- (\dtheta:\orad);
\end{scope}
}
\draw[fred] (0:\orad) -- (180:\orad);
\draw[fred] (\dtheta:\orad) -- ({10*\dtheta}:\orad);
\draw[fred] ({2*\dtheta}:\orad) -- ({6*\dtheta}:\orad);
\draw[fred] ({3*\dtheta}:\orad) -- ({9*\dtheta}:\orad);
\draw[fred] ({4*\dtheta}:\orad) -- ({13*\dtheta}:\orad);
\draw[fred] ({5*\dtheta}:\orad) -- ({11*\dtheta}:\orad);
\draw[fred] ({8*\dtheta}:\orad) -- ({12*\dtheta}:\orad);
\foreach \v in {0,1,2,...,13}
{
\begin{scope}[rotate={\v*\dtheta}]
\node[fred,fill={rgb:green,3;blue,3;white,9}] at (0:\orad) {};
\end{scope}
}
\end{scope}
\newcommand\dth{45}
\newcommand\arad{4.0}
\newcommand\brad{2.0}
\begin{scope}[xshift=10.0cm]
\begin{scope}[rotate=90]
\foreach \theta in {0,45,...,315}
{
\begin{scope}[rotate={\theta}]
\draw[fred] (0:\arad) -- (\dth:\arad);
\draw[fred] (0:\brad) -- ({3*\dth}:\brad);
\draw[fred] (0:\arad) -- (0:\brad);
\end{scope}
}
\foreach \theta in {0,45,...,315}
{
\begin{scope}[rotate={\theta}]
\node[fred,fill={rgb:green,3;blue,3;white,9}] at (0:\arad) {};
\node[fred,fill={rgb:green,3;blue,3;white,9}] at (0:\brad) {};
\end{scope}
}
\end{scope}
\end{scope}
\newcommand\dthpap{20}
\newcommand\crad{4.0}
\begin{scope}[xshift=20.0cm]
\begin{scope}[rotate=100]
\foreach \theta in {0,20,40,...,340}
{
\begin{scope}[rotate=\theta]
\draw[fred] (0:\crad) -- (\dthpap:\crad);
\end{scope}
}
\foreach \theta in {0,120,240}
{
\begin{scope}[rotate=\theta]
\draw[fred] (0:\crad) -- ({7*\dthpap}:\crad);
\draw[fred] ({2*\dthpap}:\crad) -- ({15*\dthpap}:\crad);
\draw[fred] ({10*\dthpap}:\crad) -- ({17*\dthpap}:\crad);
\end{scope}
}
\foreach \theta in {0,20,40,...,340}
{
\begin{scope}[rotate={\theta}]
\node[fred,fill={rgb:green,3;blue,3;white,9}] at (0:\crad) {};
\end{scope}
}
\end{scope}
\end{scope}
\end{tikzpicture}
\caption{Cubic graphs of orders 14, 16 and 18 having diameter $4$ and AC =
$1.2679$. The middle graph is the M\"{o}bius-Kantor graph and the graph on the
right is the Pappus graph.}%
\label{figd3D4}%
\end{figure}
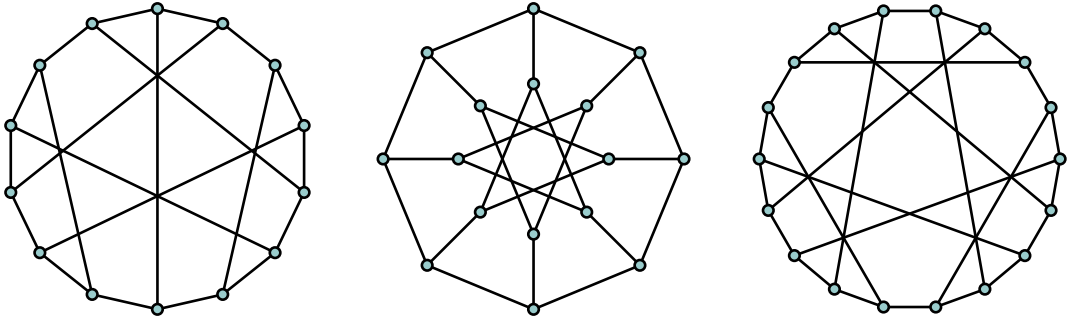

There are three graphs that achieve the AC bound of $3-\sqrt{3}\approx1.2679$.
The graphs have orders $14$, $16$ and $18$. The graph of order $14$ is the
\textquotedblleft Crossing number 3H\textquotedblright\ graph from
\cite{pegg2009crossing} and is shown in figure \ref{fig:max}(a). The graph of
order $16$ is the M\"{o}bius-Kantor graph and the graph of
order $18$ is the Pappus graph; they are shown in figure \ref{figd3D4}. 

\subsection{Degree 3, Diameter 5\label{sec:D35}}

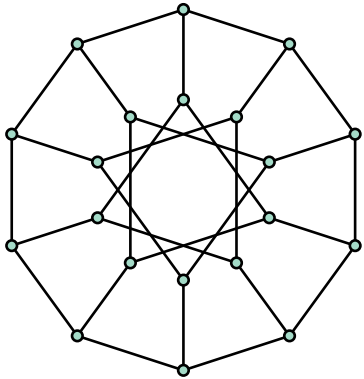
\begin{figure}[H]
\centering
\begin{tikzpicture}[scale=0.6,line width=1pt]
\newcommand\orad{4.0}
\newcommand\irad{2.0}
\newcommand\dtheta{36}
\begin{scope}[rotate=90]
\foreach \theta in {0,36,...,324}
{
\begin{scope}[rotate={\theta}]
\draw[fred] (0:\orad) -- (\dtheta:\orad);
\draw[fred] (0:\irad) -- ({3*\dtheta}:\irad);
\draw[fred] (0:\irad) -- (0:\orad);
\end{scope}
}
\foreach \theta in {0,36,...,324}
{
\begin{scope}[rotate={\theta}]
\node[fred,fill={rgb:green,3;blue,2;white,9}] at (0:\orad) {};
\node[fred,fill={rgb:green,3;blue,2;white,9}] at (0:\irad) {};
\end{scope}
}
\end{scope}
\end{tikzpicture}
\caption{The Desargues graph with diameter $5$ and AC exactly $1.0$.}%
\label{figd3D5}%
\end{figure}

By Theorem \ref{thm:Dodd}, maximal such graph must have exactly 20 vertices.
There are a total of 510489 cubic graphs on 20 vertices
\cite{coolsaet2023house, mckay1998fast}. Of these, there are exactly
\emph{five }diameter-5 graphs of that achieve the AC bound of $1$. All five have girth $6$. The most symmetric example is the Desargues graph, shown
in Figure~\ref{figd3D5}, which has vertex and edge transitive automorphism
group of order $240.$ Its cospectral mate is also one of the five, shown in
Figure \ref{fig:max}(b).

\subsection{Degree 3, Diameter 6\label{sec:D36}}

We found graphs achieving the AC bound of $3-\sqrt{5}\approx0.7639$ for all
(even) orders from $32$ to $42$, inclusive (note that the lower bound on the
order of $n=28$ from Theorem \ref{thm:Dodd} is not attained). There are two
graphs each for the six orders, except for $38$. There is only one graph
achieving the AC bound for order $38$. Graphs for orders 32 and 34 all have
girth 7. Graphs for orders 36, 38, 40 and 42 all have girth 8.

\begin{figure}[H]
\centering\includegraphics[width=0.35\textwidth]{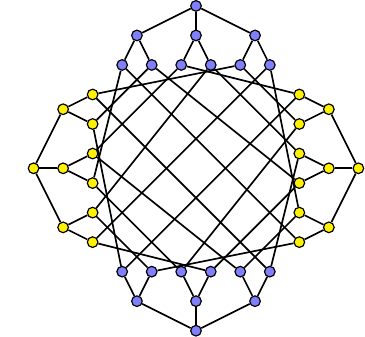} \caption{One of
the two maximal cubic graphs with $n=40$, $D=6$. It has automorphism group of
order 480.}%
\label{fig:n40}%
\end{figure}

One of the two graphs of order 40 has an automorphism group of size 480. It is
shown in figure \ref{fig:n40} .

\begin{figure}[H]
\centering
\includegraphics[width=0.7\textwidth]{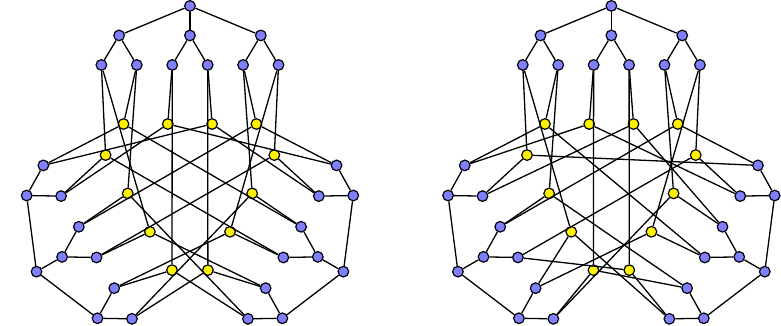} \caption{Two maximal
cospectral graphs on 42 vertices with $D=6$. Each yellow vertex is adjacent to
all three Bethe trees in blue. The graph on the left is symmetric under
rotation by 120 degrees as well as reflection in vertical axis, and has
automorphism group of order 48. The graph on the right is not, and has
automorphism group of order 24. Both have girth 8. Their spectrum is $\pm
3,\pm\sqrt{5}^{(8)},\pm2,\pm\sqrt{2}^{(4)},\pm1^{(2)},0^{(10)}$. }%
\label{figd3D6n42}%
\end{figure}

The two graphs of order 42 have the same spectrum. They both exhibit a tripe-tree
structure as shown in figure \ref{figd3D6n42}. 

\subsection{Degree 3, Diameter 7\label{sec:D37}}

\begin{figure}[H]
\centering
\begin{tikzpicture}[scale=0.5,line width=1pt]
\tikzstyle{fred}=[draw=black, fill={rgb:blue,3;green,3;white,9},
shape=circle, minimum height=1.5mm, inner sep=1, text=black]
\foreach \theta/\hue in {90/green,270/red}
{
\begin{scope}[rotate=\theta]
\coordinate (u0) at (4,0);
\coordinate (u1) at (3,-4);
\coordinate (u2) at (3,0);
\coordinate (u3) at (3,4);
\coordinate (u4) at (2,-5);
\coordinate (u5) at (2,-3);
\coordinate (u6) at (2,-1);
\coordinate (u7) at (2,1);
\coordinate (u8) at (2,3);
\coordinate (u9) at (2,5);
\coordinate (u10) at (1,-5.5);
\coordinate (u11) at (1,-4.5);
\coordinate (u12) at (1,-3.5);
\coordinate (u13) at (1,-2.5);
\coordinate (u14) at (1,-1.5);
\coordinate (u15) at (1,-0.5);
\coordinate (u16) at (1,0.5);
\coordinate (u17) at (1,1.5);
\coordinate (u18) at (1,2.5);
\coordinate (u19) at (1,3.5);
\coordinate (u20) at (1,4.5);
\coordinate (u21) at (1,5.5);
\draw[fred] (u0) -- (u1);
\draw[fred] (u0) -- (u2);
\draw[fred] (u0) -- (u3);
\draw[fred] (u1) -- (u4);
\draw[fred] (u1) -- (u5);
\draw[fred] (u2) -- (u6);
\draw[fred] (u2) -- (u7);
\draw[fred] (u3) -- (u8);
\draw[fred] (u3) -- (u9);
\draw[fred] (u4) -- (u10);
\draw[fred] (u4) -- (u11);
\draw[fred] (u5) -- (u12);
\draw[fred] (u5) -- (u13);
\draw[fred] (u6) -- (u14);
\draw[fred] (u6) -- (u15);
\draw[fred] (u7) -- (u16);
\draw[fred] (u7) -- (u17);
\draw[fred] (u8) -- (u18);
\draw[fred] (u8) -- (u19);
\draw[fred] (u9) -- (u20);
\draw[fred] (u9) -- (u21);
\node[fred] at (u0) {};
\node[fred] at (u1) {};
\node[fred] at (u2) {};
\node[fred] at (u3) {};
\node[fred] at (u4) {};
\node[fred] at (u5) {};
\node[fred] at (u6) {};
\node[fred] at (u7) {};
\node[fred] at (u8) {};
\node[fred] at (u9) {};
\node[fred,fill=\hue] at (u10) {};
\node[fred,fill=\hue] at (u11) {};
\node[fred,fill=\hue] at (u12) {};
\node[fred,fill=\hue] at (u13) {};
\node[fred,fill=\hue] at (u14) {};
\node[fred,fill=\hue] at (u15) {};
\node[fred,fill=\hue] at (u16) {};
\node[fred,fill=\hue] at (u17) {};
\node[fred,fill=\hue] at (u18) {};
\node[fred,fill=\hue] at (u19) {};
\node[fred,fill=\hue] at (u20) {};
\node[fred,fill=\hue] at (u21) {};
\end{scope}
}
\end{tikzpicture}
\caption{Starting configuration for generating graphs with diameter $7$ and
maximum AC.}%
\label{figd3D7}%
\end{figure}
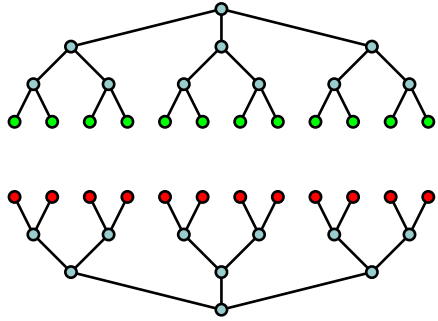

For diameter $7$ we found $45$ graphs that achieve the AC bound. All of these
graphs are bipartite with girth $8$, and were constructed using a method that
seems useful in this context. The procedure begins with two copies of the
degree $3$ Moore tree of depth $3$, as shown Figure~\ref{figd3D7}. To complete
the graph we consider only edges that join leaf vertices in of the Moore trees
to leaf vertices in the other tree. In the figure, these would be edges
joining red vertices and green vertices. There are only $144$ such edges, so
complete search for such graphs can be done quickly. Such a search results in
the $45$ graphs mentioned earlier. Among them, the highest automorphism group
has order 48, represented by a single graph.

\subsection{Degree 3, Diameter 8\label{sec:D38}}

In this case, max AC=$3-\sqrt{6}\approx0.5505.$ We found $18$ cubic graphs of
order $68$ that attain this bound. All of these graph have girth $9$. None of
the graph are particularly symmetric, with automorphism groups ranging in size
from $1$ to $24$. The graphs were constructed by searching girth $9$ graphs of
order $9$, and checking the diameter after a girth $9$ graph was constructed.

\subsection{Degree 3, Diameter 9\label{sec:D39}}

For diameter $9$, by Theorem \ref{thm:Dodd}, all maximal graphs have exactly 92 vertices. We again used the double tree method outlined in \S \ref{sec:D37}. We generated a total of around 1500 maximal graphs which took several days of computing on 6 processors at once. Of these, 481 were distinct (non-isomorphic).
All of them had girth 10.  The following table lists the statistics for group size of this collection: 
$$
\begin{tabular}
[c]{|c||c|c|c|c|c|c|c|c|c|}\hline
group size & 48 & 32 & 24 & 16 & 12 & 8 & 4 & 2 & 1\\\hline
\#graphs & 1 & 3 & 2 & 7 & 2 & 34 & 101 & 242 & 86\\\hline
\end{tabular}
$$
Note that this is not an exhautsive list, but we estimate that the actual number is close to 500.
Note that the cubic cage of girth $10$ has order $70$ \cite{cagesurvey}.

\subsection{Degree 4, Diameter 3\label{sec:D43}}

\begin{figure}[H]
\centering
\begin{tikzpicture}[scale=0.6,line width=1pt, rotate=90]
\tikzstyle{fred}=[draw=black, fill=yellow, shape=circle, minimum height=1.5mm, inner sep=1, text=black]
\newcommand\orad{4.0}
\newcommand\dtheta{36.0}
\begin{scope}[rotate=90]
\foreach \theta in {0,1,2,...,9}
{
\begin{scope}[rotate={\theta*\dtheta}]
\draw[fred] (0:\orad) -- (\dtheta:\orad);
\draw[fred] (0:\orad) -- ({3*\dtheta}:\orad);
\end{scope}
}
\foreach \v in {0,1,2,...,9}
{
\begin{scope}[rotate={\v*\dtheta}]
\node[fred,fill={rgb:blue,3;green,3;white,9}] at (0:\orad) {};
\end{scope}
}
\end{scope}
\end{tikzpicture}
\caption{A $4$-regular graph with diameter $3$ and AC = $3.0$.}%
\label{figd4D3}%
\end{figure}
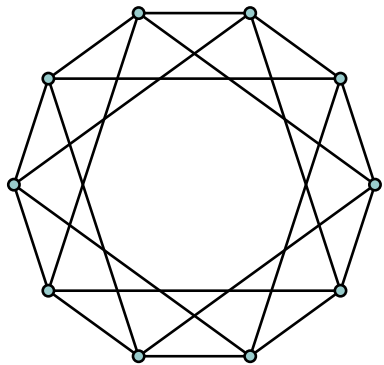

The unique maximal maximal graph with $D=3$ is the modified bipartite graph
described in \S \ref{sec:D3}. Figure \ref{figd4D3} shows its symmetric realization.

\subsection{Degree 4, Diameter 4\label{sec:D44}}

\begin{figure}[H]
\centering
\begin{tikzpicture}[scale=0.6,line width=1pt]
\tikzstyle{tree}=[fred, fill=blue!50, minimum height=2mm]
\tikzstyle{ned}=[fred,fill=yellow, minimum height=2mm]
\begin{scope}
\coordinate (u0) at (0,0);
\coordinate (u1) at (-3,2);
\coordinate (u2) at (-1,2);
\coordinate (u3) at (1,2);
\coordinate (u4) at (3,2);
\coordinate (u5) at (-5,4);
\coordinate (u6) at (-3,4);
\coordinate (u7) at (-1,4);
\coordinate (u8) at ( 1,4);
\coordinate (u9) at ( 3,4);
\coordinate (u10) at (5,4);
\coordinate (u11) at (-3,6);
\coordinate (u12) at (-1,6);
\coordinate (u13) at ( 1,6);
\coordinate (u14) at ( 3,6);
\coordinate (u15) at (0,8);
\draw[fred] (u0) -- (u1);
\draw[fred] (u0) -- (u2);
\draw[fred] (u0) -- (u3);
\draw[fred] (u0) -- (u4);
\draw[fred] (u5) -- (u1);
\draw[fred] (u5) -- (u2);
\draw[fred] (u6) -- (u1);
\draw[fred] (u6) -- (u3);
\draw[fred] (u7) -- (u1);
\draw[fred] (u7) -- (u4);
\draw[fred] (u8) -- (u2);
\draw[fred] (u8) -- (u3);
\draw[fred] (u9) -- (u2);
\draw[fred] (u9) -- (u4);
\draw[fred] (u10) -- (u3);
\draw[fred] (u10) -- (u4);
\draw[fred] (u5) -- (u11);
\draw[fred] (u5) -- (u12);
\draw[fred] (u6) -- (u11);
\draw[fred] (u6) -- (u13);
\draw[fred] (u7) -- (u11);
\draw[fred] (u7) -- (u14);
\draw[fred] (u8) -- (u12);
\draw[fred] (u8) -- (u13);
\draw[fred] (u9) -- (u12);
\draw[fred] (u9) -- (u14);
\draw[fred] (u10) -- (u13);
\draw[fred] (u10) -- (u14);
\draw[fred] (u11) -- (u15);
\draw[fred] (u12) -- (u15);
\draw[fred] (u13) -- (u15);
\draw[fred] (u14) -- (u15);
\node[tree] at (u0) {};
\node[tree] at (u1) {};
\node[tree] at (u2) {};
\node[tree] at (u3) {};
\node[tree] at (u4) {};
\node[ned] at (u5) {};
\node[ned] at (u6) {};
\node[ned] at (u7) {};
\node[ned] at (u8) {};
\node[ned] at (u9) {};
\node[ned] at (u10) {};
\node[tree] at (u11) {};
\node[tree] at (u12) {};
\node[tree] at (u13) {};
\node[tree] at (u14) {};
\node[tree] at (u15) {};
\end{scope}
\begin{scope}[xshift=12cm,yshift=4cm]
\node[ned] (1) at (-5, 0) {};
\node[ned] (2) at (-3, 0) {};
\node[ned] (3) at (-1, 0) {};
\node[ned] (4) at (1, 0) {};
\node[ned] (5) at (3, 0) {};
\node[ned] (6) at (5, 0) {};
\node[tree] (7) at (-3, 2) {};
\node[tree] (8) at (-1, 2) {};
\node[tree] (9) at (1, 2) {};
\node[tree] (10) at (3, 2) {};
\node[tree] (11) at (-3, -2) {};
\node[tree] (12) at (-1, -2) {};
\node[tree] (13) at (1, -2) {};
\node[tree] (14) at (3, -2) {};
\node[tree] (15) at (0, 4) {};
\node[tree] (16) at (0, -4) {};
\foreach \x/\y in {7/1, 8/1, 11/1, 12/1, 7/2, 9/2, 11/2, 13/2, 7/3, 10/3, 12/3, 13/3, 8/4, 9/4, 11/4, 14/4, 8/5, 10/5,
12/5, 14/5, 9/6, 10/6, 13/6, 14/6, 15/7, 15/8, 15/9, 15/10, 16/11, 16/12, 16/13, 16/14}
\path[fred] (\x) edge node {} (\y);
\end{scope}
\end{tikzpicture}%
\caption{Two maximal graphs for $d=4,D=4,n=16.$ with AC=2. Left: the Hoffman
graph. Right: the 4-Cube (Tesseract) graph}%
\label{figd4D4}%
\end{figure}
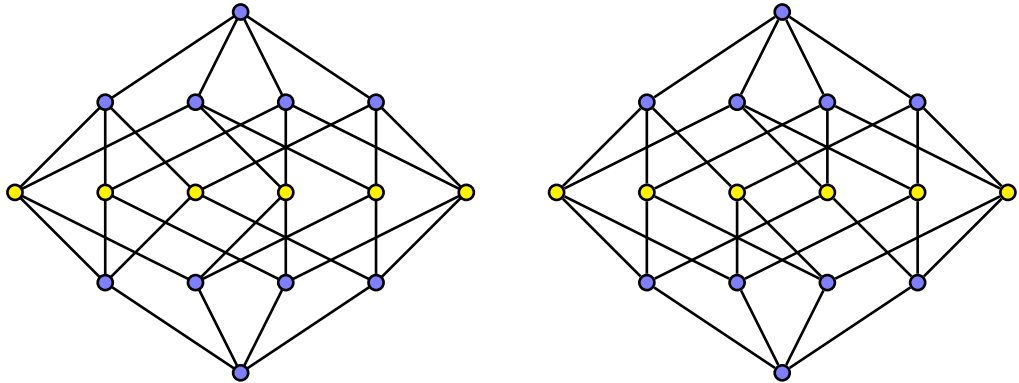

There are two maximal $4$-regular graphs of diameter $4$ on 16 vertices that
achieve the AC bound of $2.0$: the $4$-cube and its cospectral mate
\cite{hoffman}, the Hoffman graph, shown in Figure~\ref{figd4D4}. Note that
$n=16$ is the smallest possible, see Theorem \ref{thm:Deven}.

\begin{figure}[h]
\centering
\begin{tikzpicture}
\begin{scope}[scale=3]
\tikzstyle{ned}=[fred,fill=yellow!25]
\tikzstyle{ned2}=[fred,fill=blue!15]
\node[ned] (1) at (0.92388, -0.382683) {A};
\node[ned] (2) at (0.991445, -0.130526) {B};
\node[ned] (3) at (0.991445, 0.130526) {C};
\node[ned] (4) at (0.92388, 0.382683) {D};
\node[ned] (5) at (1.5, 0) {E};
\node[ned2] (6) at (0.793353, 0.608761) {F};
\node[ned2] (7) at (0.608761, 0.793353) {G};
\node[ned2] (8) at (0.382683, 0.92388) {H};
\node[ned2] (9) at (0.130526, 0.991445) {I};
\node[ned2] (10) at (0.75, 1.29904) {J};
\node[ned] (11) at (-0.130526, 0.991445) {C};
\node[ned] (12) at (-0.382683, 0.92388) {B};
\node[ned] (13) at (-0.608761, 0.793353) {D};
\node[ned] (14) at (-0.793353, 0.608761) {A};
\node[ned] (15) at (-0.75, 1.29904) {E};
\node[ned2] (16) at (-0.92388, 0.382683) {H};
\node[ned2] (17) at (-0.991445, 0.130526) {G};
\node[ned2] (18) at (-0.991445, -0.130526) {I};
\node[ned2] (19) at (-0.92388, -0.382683) {F};
\node[ned2] (20) at (-1.5, 1.83697e-16) {J};
\node[ned] (21) at (-0.793353, -0.608761) {D};
\node[ned] (22) at (-0.608761, -0.793353) {B};
\node[ned] (23) at (-0.382683, -0.92388) {A};
\node[ned] (24) at (-0.130526, -0.991445) {C};
\node[ned] (25) at (-0.75, -1.29904) {E};
\node[ned2] (26) at (0.130526, -0.991445) {I};
\node[ned2] (27) at (0.382683, -0.92388) {G};
\node[ned2] (28) at (0.608761, -0.793353) {F};
\node[ned2] (29) at (0.793353, -0.608761) {H};
\node[ned2] (30) at (0.75, -1.29904) {J};
\foreach \x/\y in {5/1, 9/1, 19/1, 27/1, 5/2, 6/2, 18/2, 29/2, 5/3, 8/3, 17/3, 28/3, 5/4, 7/4, 16/4, 26/4, 10/6, 14/6, 24/6, 10/7, 11/7, 23/7, 10/8, 13/8, 22/8, 10/9, 12/9, 21/9, 15/11, 19/11, 29/11, 15/12, 16/12, 28/12, 15/13, 18/13, 27/13, 15/14, 17/14, 26/14, 20/16, 24/16, 20/17, 21/17, 20/18, 23/18, 20/19, 22/19, 25/21, 29/21, 25/22, 26/22, 25/23, 28/23, 25/24, 27/24, 30/26, 30/27, 30/28, 30/29}
\path[ned] (\x) edge node {} (\y);
\end{scope}
\end{tikzpicture}%
\caption{Maximal graph for $d=4,D=4,n=30$ with AC=2. There are 10 triplets of
vertices, labelled A to J, that are at distance 4 from each other. Its
automorphism group is of order 720.}%
\label{fig:n30}%
\end{figure}
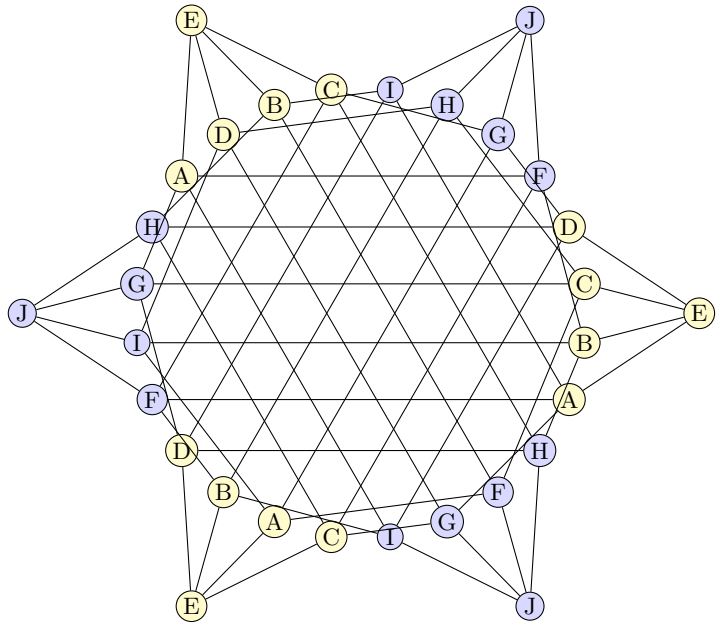

In \S \ref{sec:D4} we exhibit a maximal graph with 30 vertices based on
projective plane techniques (figure \ref{fig:n30}).

Using tables of quartic regular graphs \cite{coolsaet2023house,
meringer1999fast} we found maximal graphs up to order 32. Overall, there are
maximal graphs with 16, 22, 23, 24, 28, 30, and 32 vertices, see table in
\ref{sec:summaryD} for further information.

\subsection{Degree 4, Diameter 5\label{sec:D45}}

Twelve graphs of degree $4$ and diameter $5$ were found. All have girth $6$
and all were found using the starting configuration of two Bethe trees as in
figure \ref{figd3D7}.

\subsection{Degree 6, Diameter 4\label{sec:D64}}

In \S \ref{sec:D4} we construct a family of $D=4$ maximal graphs when $d$ is a
power of prime, having $n=2d^{2}-2$ vertices. So $d=6$ is the smallest $d$
which is not part of that family. Nonetheless, we found a maximal graph
using stochastic search algorithm.

\subsection{Summary of diameter-maximal graphs \label{sec:summaryD}}

The table below summarizes our findings. For the column \textquotedblleft%
\#graphs\textquotedblright, a number in bold indicates that the corresponding
class has been searched exhaustively, and no other graphs for the
corresponding $n$ are expected. Otherwise, it is the number we managed to
find, but there may be more. Note that we restricted the search to the girth
as specified in the table. Where possible (for $n\leq20$), we also confirmed
using exhaustive search that no maximal graphs exist for girths smaller than indicated.

\begin{center}%
\begin{tabular}
[c]{|l|l|l|l|l|l|}\hline
$d$ & $D$ & AC & $n$ & \#graphs & Comments\\\hline \hline
3 & 3 & 2 & 8 & \textbf{1} & The $3$-cube, see \S \ref{sec:D3}.\\\hline
3 & 4 & 1.2679 &
\begin{tabular}
[c]{l}%
14\\
16\\
18
\end{tabular}
&
\begin{tabular}
[c]{l}%
\textbf{1}\\
\textbf{1}\\
\textbf{1}%
\end{tabular}
&
\begin{tabular}
[c]{l}%
Graph 3H (see Figure \ref{fig:max}(a))\\
M\"{o}bius Kantor Graph\\
Pappus graph
\end{tabular}
\\\hline
3 & 5 & 1 & 20 & \textbf{5} & All have girth $6$; includes the Desargues
graph.\\\hline
3 & 6 & 0.7639 &
\begin{tabular}
[c]{l}%
32\\
34\\
36\\
38\\
40\\
42
\end{tabular}
&
\begin{tabular}
[c]{l}%
\textbf{2}\\
\textbf{2}\\
\textbf{2}\\
\textbf{1}\\
\textbf{2}\\
\textbf{2}%
\end{tabular}
&
\begin{tabular}
[c]{l}%
Both have girth 7\\
Both have girth 7, cospectral\\
Both have girth 8\\
Girth 8\\
Both have girth 8\\
Both have girth 8, cospectral
\end{tabular}
\\\hline
3 & 7 & 0.6571 & $44$ & \textbf{45} & All of have girth 8.\\\hline
3 & 8 & 0.5505 &
\begin{tabular}
[c]{l}%
68\\
80\\
90
\end{tabular}
&
\begin{tabular}
[c]{l}%
12\\
1\\
3
\end{tabular}
&
\begin{tabular}
[c]{l}%
Two of girth 8 and ten of girth 9\\
Girth 10\\
Girth 10
\end{tabular}
\\\hline
3 & 9 & 0.4965 & 92 & $481$ & All have girth 10\\\hline\hline
4 & 3 & 3 & 10 & \textbf{1} & See \S \ref{sec:girth34}.\\\hline
4 & 4 & 2 &
\begin{tabular}
[c]{l}%
16\\
17\\
19-21\\
22\\
23\\
24\\
28\\
30\\
32
\end{tabular}
&
\begin{tabular}
[c]{l}%
\textbf{6}\\
\textbf{0}\\
\textbf{0}\\
\textbf{3}\\
\textbf{2}\\
\textbf{2}\\
\textbf{1}\\
\textbf{1}\\
\textbf{1}%
\end{tabular}
&
\begin{tabular}
[c]{l}%
Girth 4; group sizes:\ 6,8,12,32,48 (Hoffman), 384 (Tesseract)\\
Girth 4:\ 193900 graphs, none maximal\\
Girth 5\\
Girth 5, group sizes 2, 4, 8\\
Girth 5; group sizes 1, 4\\
Girth 5; co-spectral, group sizes 16, 16\\
Unique graph of girth 6 on 28 vertices\\
Girth 6, see \S \ref{sec:D4} and figure \ref{fig:n30}\\
Girth 6
\end{tabular}
\\\hline\hline
6 & 4 & 3.5505 & 44 & $\geq1$ & Not part of $D=4$ family of \S \ref{sec:D4}%
\\\hline
\end{tabular}

\end{center}

\section{Maximum AC for given girth and order\label{sec:girth}}

As mentioned in Theorem \ref{thm:girth}, girth-maximal graphs must necessarily
be Moore graphs. This gives a severe restriction on existence of girth-maximal
graphs. Conversely, all the Moore graphs we considered appear to be
girth-maximal. This includes girth-6 projective plane family (see
\S \ref{sec:girth6}), Peterson graph on 10 vertices of girth 5, 
the cubic Tutte cage on 30 vertices of
girth 8, and the cubic Benson cage of girth 12 on 126 vertices \cite{cagesurvey}. It
remains an open question as to whether all known Moore graphs are girth-maximal.

We used the data from House of Graphs website \cite{coolsaet2023house} which
contains complete enumeration of cubic graphs up to orders 64, particularly
for higher girths. For each combination of order and girth, we computed the
maximum AC. This is recorded in the tables below.%

\begin{tabular}
[c]{|c|c|c|c|}\hline
$g$ & $n$ & max AC & Comments\\\hline\hline
3 &
\begin{tabular}
[c]{c}%
4\\
6\\
8\\
10\\
12\\
14\\
16\\
18\\
\end{tabular}
&
\begin{tabular}
[c]{c}%
4\\
2\\
1.438447\\
1.120614\\
1\\
0.885092\\
0.82259\\
0.763932\\
\end{tabular}
&
\begin{tabular}
[c]{c}%
Complete graph, maximal\\
\\
\\
\\
\\
\\
\\
\\
\end{tabular}
\\\hline
4 &
\begin{tabular}
[c]{c}%
6\\
8\\
10\\
12\\
14\\
16\\
18\\
20\\
\end{tabular}
&
\begin{tabular}
[c]{c}%
3\\
2\\
1.438447\\
1.267949\\
1.068732\\
1\\
0.903097\\
0.845793\\
\end{tabular}
&
\begin{tabular}
[c]{l}%
Modified bipartite, maximal, see \S \ref{sec:girth34}\\
\\
\\
\\
\\
\\
\\
\\
\end{tabular}
\\\hline
5 &
\begin{tabular}
[c]{c}%
10\\
12\\
14\\
16\\
18\\
20\\
22
\end{tabular}
&
\begin{tabular}
[c]{c}%
2\\
1.467911\\
1.289171\\
1.172909\\
1.043705\\
1\\
0.913969
\end{tabular}
&
\begin{tabular}
[c]{l}%
Petersen graph, maximal\\
\\
\\
\\
\\
\\
\\
\end{tabular}
\\\hline
\end{tabular}
\ \
\begin{tabular}
[c]{|c|c|c|c|}\hline
$g$ & $n$ & {max AC} & {Comments}\\\hline\hline
6 &
\begin{tabular}
[c]{c}%
14\\
16\\
18\\
20\\
22\\
24\\
\end{tabular}
&
\begin{tabular}
[c]{c}%
1.585786\\
1.267949\\
1.267949\\
1.064568\\
1\\
1\\
\end{tabular}
&
\begin{tabular}
[c]{c}%
PG(2,2), maximal, see \S \ref{sec:girth6}\\
\\
\\
\\
\\
\\
\end{tabular}
\\\hline
7 &
\begin{tabular}
[c]{c}%
24\\
26\\
28\\
30\\
32\\
30
\end{tabular}
&
\begin{tabular}
[c]{c}%
1\\
0.94737\\
1\\
0.844084\\
0.864221\\
1
\end{tabular}
&
\begin{tabular}
[c]{c}%
\\
\\
\\
\\
\\
\end{tabular}
\\\hline
8 &
\begin{tabular}
[c]{c}%
30\\
34\\
36\\
38\\
40\\
42
\end{tabular}
&
\begin{tabular}
[c]{c}%
1\\
{0.78568}\\
{0.763932}\\
{0.763932}\\
{0.763932}\\
{0.763932}%
\end{tabular}
&
\begin{tabular}
[c]{c}%
Tutte graph, maximal\\
\\
\\
\\
\\
$\ \ $%
\end{tabular}
\\\hline
9 &
\begin{tabular}
[c]{c}%
58\\
60\\
62\\
64
\end{tabular}
&
\begin{tabular}
[c]{c}%
0.63766\\
0.697224\\
0.603671\\
0.633832
\end{tabular}
&
\begin{tabular}
[c]{l}%
\\
\\
\\
Next 3: 0.6031,0.6025,0.6011
\end{tabular}
\\\hline
12 & 126 & {0.5505} & {Benson graph, maximal}\\\hline
\end{tabular}

\begin{center}

\end{center}

\section{Maximal graph families}\label{sec:families}

Here, we exhibit families of maximal graphs for infinitely many $d$: girth
3,4, and 6, and diameters 3 and 4.

\subsection{Girth 3 and 4\label{sec:girth34}}

The trivial cases are the girth-3 and girth-4 maximal graphs which are
complete graph on $d+1$ vertices, and complete-bipartite graphs on $2d$
vertices. Their AC is well known to be $d+1$ and $d$ respectively, which are maximal.

\subsection{Diameter 3\label{sec:D3}}

This is a modified bipartite graph on $2d+2$ vertices, illustrated here for
the case of $d=4:$%
\[
\includegraphics[width=0.2\textwidth]{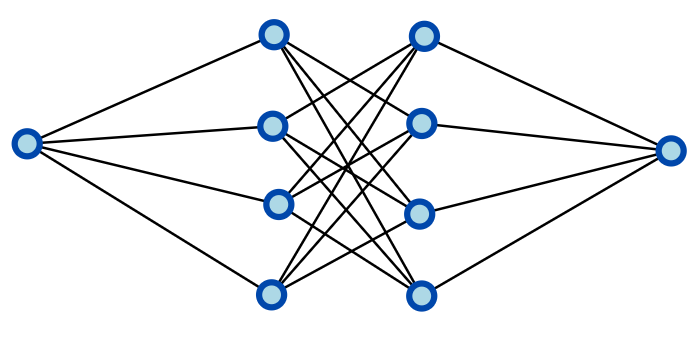}
\]
It is constructed by starting with complete bipartite graph on $2d$ vertices,
removing $d$ matching edges, and adding two extra vertices, each with $d$
edges connecting it to $d$ vertices of each of component. This graph has
spectrum $\pm d$ (once)$,\ \pm1$ ($d$ times each), so that $AC=d-1,$ which is
maximal for $D=3.$

This maximal graph is in fact unique, which can be seen as follows. Maximal
graph of odd diameter $D=2K-1$ must consist of two Bethe trees of $K$ levels
whose leafs are connected to each-other. In the case of $D=3$, the two trees
are simply two star graphs with $d$ leafs. Each leaf of the left star is
connected to $d-1$ leafs of the right star and vice-versa. In other words, for
each leaf on the left star, there is exactly one leaf on the right star to
which it is not connected and vice-versa. Moreover, this is a one-to-one
correspondence (otherwise there would be more than one leaf on the left star
not connected to a single leaf on the right star). This shows that all maximal
graphs with $D=3$ must be isomorphic.

\subsection{Girth 6.\label{sec:girth6}}

It is well known that a projective plane\textbf{ }$PG(2,q)$, $q$ a prime
power, is in fact a degree $d=q+1$ regular graph of girth 6, and is indeed one
of the few Moore graphs \cite{cagesurvey}. Its spectrum is easy to compute,
see for example \cite{ryser1955geometries, hoffman1965line,
godsil1977spectrum}. From there it follows that $AC=d-\sqrt{d-1}$ which is
maximal for girth 6. For completeness, we reproduce these arguments here.

Let us recall the construction of a projective plane graph. Given a degree
$d,$ consider a finite field $F$ of size $q=d-1$ (such a field only exists
when $q$ is a prime power). A \emph{line }$L\in PG(2,q)$ is a non-zero tuple
$\left(  a,b,c\right)  $ with $a,b,c\in F$, modded by an equivalence relation
corresponding to scalar multiplication. More explicitly, a line can be
represented uniquely by rescaling the first non-zero coordinate to
one:\ either $\left(  1,b,c\right)  $ with $b,c\in F$ or $\left(
0,1,c\right)  $ with $c\in F$, or $\left(  0,0,1\right)  .$ Correspondingly,
there are $q^{2}+q+1$ such lines. A\ point $P=(x,y,z)$ in $PG(2,q)$ has the
same form as a line. Then $PG(2,q)$ is a bipartite graph having $2(q^{2}+q+1)$
vertices. Half of the vertices are lines $L$, half are points $P$, and there
is an edge between $L$ and $P$ if and only if $L\cdot P\equiv0.$

One can easily check that $PG(2,q)$ is regular of degree $d=q+1;$ has
$n=2d^{2}-2d+2$ vertices, has no four-cycles (and therefore has girth $g=6),$
and has diameter $D=3$. Let us compute its spectrum following
\cite{ryser1955geometries, hoffman1965line, godsil1977spectrum}.

The adjacency matrix has the form $A=\left[
\begin{array}
[c]{cc}%
0 & B\\
B & 0
\end{array}
\right]  $ where $B_{ik}=\left\{
\begin{array}
[c]{c}%
1\text{, if }L_{i}\perp P_{k}\\
0\text{, otherwise}%
\end{array}
\right.  .$ Note that $B$ is symmetric since lines and points are identical
and interchangeable in this geometry. Correspondingly, the eigenvalues of $A$
are given by $\lambda=\pm\sqrt{\mu}$, where $\mu$ is an eigenvalue of
$M=B^{2}.$ Note that $M_{ij}=\sum_{k}B_{ik}B_{kj}$ so that $M_{ij}$ is number
of points $P$ which are simultaneously orthogonal to both lines $L_{i}$ and
$L_{j}$. It is easy to check that%
\[
M_{ij}=\left\{
\begin{array}
[c]{c}%
q+1,\text{ if }i=j\\
1\text{ otherwise}%
\end{array}
\right.
\]
Correspondingly, the eigenvalues of $M$ are $\mu=q=d-1$ (with multiplicity
$q^{2}+q$), and $q^{2}+2q+1=d^{2}$ with multiplicity one. So the spectrum of
$A$ consists of four eigenvalues: $\pm d,$ and $\pm\sqrt{d-1}$ with
multiplicity $d(d-1).$ It follows that $AC=d-\sqrt{d-1}.$

\subsection{Diameter 4.\label{sec:D4}}

Here, we will construct a $d-$ regular graph $G$ which is a subgraph of
$PG(2,d)$ (with $d$ a prime power). Its order is $2d^{2}-2.$ This graph is
likely to be the same as girth-6 graph of same order from
\cite{abreu2006minimal, araujo2010finding, araujo2011constructions}, although
we use a different construction here to compute its spectrum and girth.

Consider the subset of lines and points of $PG(2,d)$ of the form $\left(
1,b,c\right)  $, where one of $b,c$ are non-zero. For example when $d=3,$
there are 8 such lines and points, namely:%
\begin{equation}
\left(  1,0,1\right)
,\ (1,0,2);\ \ \ (1,1,0),\ (1,2,0);\ \ \ (1,1,1),\ (1,2,2);\ \ \ (1,1,2),\ (1,2,1).
\label{223}%
\end{equation}
It is easy to see that such graph is regular of degree $d;$ has $n=2d^{2}-2$
vertices, its girth is $g=6$ and its diameter is $D=4.$ We start by showing
the latter here.

Consider two distinct lines $L_{1}=\left(  1,a_{1},b_{1}\right)  $ and
$L_{2}=\left(  1,a_{2},b_{2}\right)  .$ They are adjacent to the same point
$P=\left(  1,x,y\right)  $ if and only if $\left(
\begin{array}
[c]{cc}%
a_{1} & b_{1}\\
a_{2} & b_{2}%
\end{array}
\right)  \left(
\begin{array}
[c]{c}%
x\\
y
\end{array}
\right)  =\left(
\begin{array}
[c]{c}%
-1\\
-1
\end{array}
\right)  .$ If this system has a solution, then the distance between these two
lines is 2. In the opposite case, we have that $\left(  a_{2},b_{2}\right)
=c\left(  a_{1},b_{1}\right)  $ for some $c\in F.$ In this case, pick a point
$P$ perpendicular to $L_{1}.$ This point has $d-1$ other lines that are
perpendicular to it. Pick one such line, call it $L_{3}=\left(  1,a_{3}%
,b_{3}\right)  .$ Note that $\left(  a_{3},b_{3}\right)  \neq c\left(
a_{2},b_{2}\right)  $. for any $c\in F.$ But then $\operatorname{dist}%
(L_{2},L_{3})=2=\operatorname{dist}(L_{3},L_{1})$ so that $\operatorname{dist}%
(L_{1},L_{2})=4.$ Similar argument shows that $\operatorname{dist}(L,P)\leq3$
for any line $L$ and point $P.$

Next we compute the spectrum of $G.$ As before, its spectrum is given by
$\lambda=\pm\sqrt{\mu}$, where $\mu$ is an eigenvalue of matrix $M=B^{2};$
where $B_{ij}=\left\{
\begin{array}
[c]{c}%
1\text{, if }L_{i}\perp P_{k}\\
0\text{, otherwise}%
\end{array}
,\right.  $ with $M_{ij}$ being the number of points $P$ which are orthogonal
to both lines $L_{i}$ and $L_{j}$. For example, in the case of $d=3$ and with
lines $L_{1}\ldots L_{8}$ and points $P_{1}\ldots P_{8}$ ordered as in
(\ref{223}), the corresponding matrices are%
\[
B=\left[
\begin{array}
[c]{cccccccc}%
0 & 1 & 0 & 0 & 0 & 1 & 1 & 0\\
1 & 0 & 0 & 0 & 1 & 0 & 0 & 1\\
0 & 0 & 0 & 1 & 0 & 1 & 0 & 1\\
0 & 0 & 1 & 0 & 1 & 0 & 1 & 0\\
0 & 1 & 0 & 1 & 1 & 0 & 0 & 0\\
1 & 0 & 1 & 0 & 0 & 1 & 0 & 0\\
1 & 0 & 0 & 1 & 0 & 0 & 1 & 0\\
0 & 1 & 1 & 0 & 0 & 0 & 0 & 1
\end{array}
\right]  ;\ \ \ M=\left[
\begin{array}
[c]{cccccccc}%
3 & 0 & 1 & 1 & 1 & 1 & 1 & 1\\
0 & 3 & 1 & 1 & 1 & 1 & 1 & 1\\
1 & 1 & 3 & 0 & 1 & 1 & 1 & 1\\
1 & 1 & 0 & 3 & 1 & 1 & 1 & 1\\
1 & 1 & 1 & 1 & 3 & 0 & 1 & 1\\
1 & 1 & 1 & 1 & 0 & 3 & 1 & 1\\
1 & 1 & 1 & 1 & 1 & 1 & 3 & 0\\
1 & 1 & 1 & 1 & 1 & 1 & 0 & 3
\end{array}
\right]  .
\]
Zeros in $M$ correspond to lines that are at distance 4 from each other and
ones to lines at distance 2 from each other. To see this more generally, group
$d^{2}-1$ lines into distinct classes $C_{1},\ldots,C_{d+1}$ of $d-1$ members
each, such that lines $L_{i}=\left(  1,a_{1},b_{1}\right)  $ and
$L_{j}=\left(  1,a_{2},b_{2}\right)  $ are in the same class if and only if
$\left(  a_{2},b_{2}\right)  =c\left(  a_{1},b_{1}\right)  $ for some $c\in
F\backslash\left\{  0\right\}  .$ Then $M_{ij}=0$ iff $i\neq j$ and $i,j\in
C_{k}$ for some $k;$ otherwise $M_{ij}=1$ if $i\neq j$ and $M_{ij}=d$ if
$i=j.$

Define an eigenvector $v$ such that $\sum_{j\in C_{k}}v_{j}=0$ for each class
$C_{k},k=1\ldots d+1.$ This gives $d+1$ linear equations, so there is a total
of $\left(  d^{2}-1\right)  -(d+1)=d^{2}-d-2$ independent such eigenvectors;
and moreover it is easy to check that\ $Mv=dv.$ Thus, $M$ has an eigenvalue
$\mu=d$ of multiplicity $d^{2}-d-2.$ Next, for each $k=1\ldots d+1$, define
$v$ such that $v_{j}=V_{k}$ if $j\in C_{k};$ and moreover choose $V_{k}$ such
that $\sum_{1}^{d+1}V_{k}=0.$ Note that there are $d$ independent such
choices. Then $Mv=(d-(d-1))v=v.$ This gives an eigenvalue of $\mu=1$ of
multiplicity $d.$ Finally, the eigenvector $v=(1\ldots1)$ yields an eigenvalue
of $d^{2}.$

In conclusion, the spectrum of $G$ is $\pm d$ (once)$,\pm\sqrt{d}$
($d^{2}-d-2$ times), and $\pm1$ ($d$ times). Subsequently, $AC=d-\sqrt{d}.$

\textbf{Summary. }The following table summarizes some facts about attainable
bounds with respect to girth and diameter. Bold font indicates whether it is
maximal with respect to girth or diameter.

\begin{center}%
\begin{tabular}
[c]{|l|l|l|l|l|l|}%
\hline
$d$ & $D$ & $g$ & $n$ & AC & Comments\\\hline\hline
any $d$ & \textbf{1} & \textbf{3} & $d+1$ & $d+1$ & complete graph on $d$
vertices\\\hline
any $d$ & \textbf{2} & \textbf{4} & $2d$ & $d$ & complete-bipartite
graph\\\hline
any $d$ & \textbf{3} & 4 & $2d+2$ & $d-1$ & Modified bipartite graph\\\hline
$p^{\alpha}+1$, $p$ prime & 3 & \textbf{6} & $2d^{2}-2d+2$ & $d-\sqrt{d-1}$ &
Projective plane $PG\left(  2,d-1\right)  $\\\hline
$p^{\alpha}$, $p$ prime & \textbf{4} & 6 & $2d^{2}-2$ & $d-\sqrt{d}$ & Subset
$PG\left(  2,d\right)  ^{-}\subset PG\left(  2,d\right)  $%
\\
\hline
\end{tabular}

\end{center}

\section{Discussion and open questions}\label{sec:discuss}

We exhibited tight upper bounds for AC for regular graphs with respect to
girth or diameter. While the girth bound is attainable only by Moore graphs --
which imposes a severe restriction on $g$ -- the diameter bound is less
restricting but is nonetheless is very rarely attained. Using a combination of
stochastic algorithms and exhaustive search, we produced examples of maximal
graphs for $d=3$ and $D\leq9,$ as well as $d=4$ and $D\leq6.$ There are many
interesting open questions -- both computational and theoretical -- that we
hope the reader will be tempted to explore.

Complete lists of 
cubic and quartic graphs for small orders suggests that for 
a fixed $n$, the graph with a
maximum AC also have the maximum possible girth; see table in
\S\ref{sec:girth}. This was the key
insight that allowed us to find diameter-maximal graphs:\ we simply searched
for graphs of highest possible girth for a given order, and generated as many
such graphs as we could; then hope that a small subset of these would end up
being diameter-maximal. This leads to our first conjecture.

\begin{conj}
\label{conj:ng}
For a fixed degree $d$ and order $n$, the graph with maximum possible AC has
the maximum possible girth.
\end{conj}

In fact, the table in \S \ref{sec:girth} suggests an even stronger conjecture:

\begin{conj}
Let\textbf{ }$f(g;d,n)$ be the maximum possible AC among the graphs of given
degree $d,$ order $n$ and girth $g.$ Then $f(g;d,n)$ is an increasing function
of $g.$
\end{conj}

Our search for diameter-maximal also suggests the following conjecture:

\begin{conj}
\label{conj:Dg}
$\ \ $

\begin{itemize}
\item A diameter-maximal graph of odd diameter $D$ must have girth $g=D+1$.

\item A diameter-maximal graph of even diameter $D$ must have girth of either
$g=D,D+1$ or $D+2.$
\end{itemize}
\end{conj}

The complexity of finding maximal graphs increases tremendously with larger
$D.$ We spent significant time searching for $D=10$ maximal cubic graphs, but
did not find any as of this writing. We state this as an open problem.

\begin{openq}
Find a $D=10$ maximal cubic graph, or show it doesn't exist.
\end{openq}

We showed that girth-maximal graphs are necessarily Moore graphs, which imposes
a severe restriction on $g$.  The only possible Moore graphs
are $d=3$, $7$, and possibly $57$, and $g=5$, or else $d-1$ is a prime power and $g=6$, $8$ or $12$.
\cite{cagesurvey}). We also showed that Moore graphs with $g=6$ are maximal
for any $d$ which is a prime power plus one. In addition, we verified that
Moore graphs with $d=3$ and $g=5,8\,$and $12$ are also maximal. What aboutother $d?$

\begin{openq}
Are all Moore graphs girth-maximal?
\end{openq}

In contrast to girth-maximal, we found that diameter-maximal graphs exist for
all $D<10$ when $d=3$. Do they exist for all $D?$

\begin{openq}
Do diameter-maximal graphs exist for any $D?$
\end{openq}

We described a family of diameter-maximal graphs for $D=4$ when $d$ is a prime
power. A computer search also revealed a maximal graph with $d=6,D=4.$ The
smallest unsettled case with $D=4$ is therefore $d=10.$

\begin{openq}
Find maximal graphs with $D=4$ when $d$ is not a prime power. Find a general
family of maximal graphs for $D=5$ and higher.
\end{openq}

Finally, a big difference between odd and even diameters is that maximal
graphs for odd diameter exist only for a specific value of $n$ given in
Theorem \ref{thm:Dodd}, whereas even-diameter graphs exist for a range of
values of $n$. Theorem \ref{thm:Deven} gives the lower bound for such $n,$
although it is not always attained. What about the upper bound?

\begin{openq}
For even $D$, what is the largest $n$ that admits a maximal graph?
\end{openq}

For example, Theorem \ref{thm:Dodd} gives a lower bound of $n=28$ when
$d=3,\ D=8;$ we found examples of maximal graphs with $n=68,80,$ and $90$.
Maximal graphs for $D=9$ require $n=92.$ Do does there exist a maximal graphs
with $D=8$ and $n=92?$ Do maximal graphs exist for values of $n<68?$ For other
values of $n?$

We mostly concentrated on diameter-maximal graphs. While girth-maximal graphs
do not generally exist due to the Moore graph constraint, there is \emph{some}
graph that maximizes AC. Section \ref{sec:girth} lists some of these records
based on complete enumeration of cubic graphs of high girths
\cite{coolsaet2023house, mckay1998fast}. What about higher $n$ or $g$, where
complete enumeration is impossible?

\begin{openq}
For a given girth or given graph order $n$, find an efficient algorithm to
produce $d-$regular graph with as high AC as possible.
\end{openq}

A well-known Ghosh-Boyd algorithm \cite{ghosh2006growing} generating
well-connected graphs -- when modified to produce regular graphs
only\footnote{The original version of the Ghosh-Boyd algorithm does not
constraint graph degrees and typically results in an irregular graph
\cite{kolokolnikov2015maximizing}. We used a modified version which disallows
having degrees more than $d$ and produces a regular graph.} -- does much worse
than the optimal. For example, running modified Ghosh-Boyd algorithm 100 times
with $d=3,n=64$ yields an average AC of\ $0.45$ (std=0.018, max=0.489 over 100
simulations). Moreover almost all runs produce girth 5 or 6; none produced
more than 7. By contrast, the record graph (see \S \ref{sec:girth})\ is a
graph of girth 9 and AC=0.6338; the average AC of all 1408813 girth-9 such
graphs is $0.53$ (std=0.018), quite far from the Ghosh-Boyd result. We also
generated random cubic graphs on 64 vertices. These have average AC of 0.25
(std=0.047, max=0.367 over 100 simulations). Finding the record graph is an
extremely time-consuming task; it takes thousands of computer hours to find the record
AC using the full enumeration compared to a fraction of a second
to run the Ghosh-Boyd algorithm. Finding an \emph{efficient} algorithm that
does significantly better than Ghosh-Boyd is an important problem.

\appendix

\section{Code for maximal graph generation}

Many of the graph discussed in the paper were found by a search procedure
outlined below. The C code is available for download \cite{code}, along with a
collection of maximal graphs.

The procedures uses a $Graph$ data structure that contains

\begin{itemize}
\item The adjacency list,

\item The adjacency matrix,

\item The edge list,

\item The edge count, denoted $edgecount$,

\item The target, i.e., the number of edges in a completed graph.
\end{itemize}

The procedure begins by initializing the graph. This initialization either
results in a graph with no edges or else a forest. Examples of forests that
might be used include a Moore tree or pair of trees like the ones shown in
Figure~\ref{figd3D7}.

The main loop of the procedure begins by calling the \texttt{Makelist}
procedure which creates a list of edges that could be added to the graph
without violating the degree and girth constraints. If the initial graph has
no edges, then for the first iteration this list will contain all possible
edges. The list is then sorted by the sum of the degrees of the two vertices
in the edge, in decreasing order. Candidate edges with the same degree sum are
sorted randomly. This sorting by decreasing degree sum is the single most
important idea in the procedure. Without it, the procedure will be successful
only on very small graphs. With it, new results can be obtained. For example,
the $(4,7)$-cage was found using this method \cite{exoo2011computational}.

At the end of the for-loop in the procedure, we have either completed the
graph, in which case we are done, or we remove a small number of edges and try
again. The number of edges removed is typically very small (usually $1$), but
if no increase in the maximum number of edges attained is realized for a large
number of iteration through the while loop, the number of edges removed can be
slowly increased, until reaching some maximum value, after which it is reset
to the small value.

\begin{algorithm}[H]
\begin{small}
\caption{Graph Search}
\label{graphalgo}
\begin{algorithmic}
\algnewcommand{\BigComment}[1]{\State \(\triangleright\) #1}
\Procedure{Search}{$graph$,$initialState$}
\State Initialize Graph
\While{$edgecount < target$} \Comment{The $target$ is $degree*order/2$}
\State $feasible \gets Makelist(graph)$
\State $SortEdges(feasible,graph)$
\For{each edge $e$ in $feasible$}
\State SortEdges($graph$)
\If{EdgeCheck($e$,$graph$)}
\State AddEdge($e$, $graph$)
\EndIf
\State{Remove $e$ from $feasible$}
\EndFor
\If{$edgecount == target$}
\State{Return}
\EndIf
\State{Remove a few random edges}
\EndWhile
\EndProcedure
\end{algorithmic}
\end{small}
\end{algorithm}

\bibliographystyle{pnas2009}
\bibliography{graph}

\end{document}